\documentclass[11pt,letterpaper]{article}
\usepackage[utf8]{inputenc}
\usepackage[T1]{fontenc}
\usepackage{lmodern}
\usepackage{natbib}
\usepackage{bibunits}
\defaultbibliographystyle{apalike}
\setcitestyle{authoryear,open={(},close={)},citesep={,}}
\usepackage{amsmath,amssymb,amsthm,mathtools}
\usepackage{aliascnt}
\usepackage{pgfplots}
\usepgfplotslibrary{groupplots}
\pgfplotsset{compat=1.18}
\definecolor{paperGlobal}{RGB}{38,77,126}
\definecolor{paperLocal}{RGB}{183,102,39}
\definecolor{paperEL}{RGB}{78,121,105}
\newcommand{\paperfigurewidth}{\ifdefined\figwidth\figwidth=\linewidth\fi}
\pgfplotsset{
 paperaxis/.style={
  scale only axis,width=3.9cm,height=2.0cm,
  font=\fontsize{8}{9}\selectfont,
  title style={font=\fontsize{8}{9}\selectfont},
  label style={font=\fontsize{8}{9}\selectfont},
  tick label style={font=\fontsize{7.5}{8.5}\selectfont},
  xmode=log,log basis x=2,xmin=40,xmax=250,
  xtick={50,100,200},xticklabels={50,100,200},
  ymin=86,ymax=101,ytick={90,95,100},
  axis x line*=bottom,axis y line*=left,
  tick align=outside,tick style={black!65,thin},
  axis line style={black!65,thin},grid=none,clip=true,
  error bars/y dir=both,error bars/y explicit,
  error bars/error bar style={line width=.45pt},
  error bars/error mark options={mark size=1.4pt,line width=.45pt}
 },
 paperg/.style={only marks,mark=*,mark size=1.6pt,
  color=paperGlobal,line width=.5pt},
 paperl/.style={only marks,mark=square*,mark size=1.6pt,
  color=paperLocal,line width=.5pt},
 paperel/.style={only marks,mark=diamond,mark size=2pt,
  color=paperEL,line width=.5pt},
 paperwald/.style={only marks,mark=triangle,mark size=2pt,
  color=black!75,line width=.5pt}
}

\usepackage{dsfont}
\usepackage{float}
\usepackage{flafter}
\usepackage[dvipsnames]{xcolor}
\usepackage[letterpaper,margin=1in]{geometry}
\usepackage{setspace}
\usepackage[shortlabels]{enumitem}
\setlist{topsep=0.35em,itemsep=0.15em,parsep=0pt,partopsep=0pt}
\usepackage[expansion=false]{microtype}
\usepackage{titling}
\usepackage{abstract}
\usepackage[colorlinks=true,linkcolor=blue,citecolor=blue,urlcolor=blue,hypertexnames=false]{hyperref}
\usepackage[nameinlink,noabbrev]{cleveref}
\newcommand{\maincref}[1]{\cref{#1}}
\newcommand{\mainCref}[1]{\Cref{#1}}
\newcommand{\mainequationref}[1]{\eqref{#1}}
\newcommand{\suppcref}[1]{\cref{#1}}
\newcommand{\suppCref}[1]{\Cref{#1}}
\pretitle{\begin{center}\huge\bfseries}
\posttitle{\par\end{center}\vskip 0.5em}
\preauthor{\begin{center}\Large\scshape}
\postauthor{\par\end{center}\vspace{-0.6em}}
\predate{\begin{center}\normalsize\itshape}
\postdate{\par\end{center}\vspace{-0.45em}}
\renewcommand{\and}{\par}

\newcommand{\paperkeywords}[1]{\par\noindent\textbf{Keywords:} #1\par}
\theoremstyle{plain}
\newtheorem{theorem}{Theorem}[section]
\newaliascnt{lemma}{theorem}
\newtheorem{lemma}[lemma]{Lemma}
\aliascntresetthe{lemma}
\newaliascnt{corollary}{theorem}
\newtheorem{corollary}[corollary]{Corollary}
\aliascntresetthe{corollary}
\newaliascnt{proposition}{theorem}
\newtheorem{proposition}[proposition]{Proposition}
\aliascntresetthe{proposition}
\theoremstyle{definition}
\newaliascnt{definition}{theorem}

\aliascntresetthe{definition}
\newtheorem{assumption}{Assumption}
\newtheorem{condition}{Condition}[section]
\theoremstyle{remark}
\newtheorem{remark}{Remark}[section]

\providecommand{\pr}{\mathop{\mathbb P}\nolimits}
\providecommand{\E}{\mathop{\mathbb E}\nolimits}
\providecommand{\Var}{\operatorname{Var}}
\providecommand{\Cov}{\operatorname{Cov}}
\providecommand{\argmin}{\operatorname*{arg\,min}}
\providecommand{\wh}{\widehat}
\providecommand{\cQ}{\mathcal Q}

\hypersetup{allcolors=blue}
\DeclareRobustCommand{\eqrefs}[1]{\labelcref{#1}}
\DeclareRobustCommand{\crefplainname}[1]{\textcolor{black}{#1}}
\crefname{theorem}{\crefplainname{Theorem}}{\crefplainname{Theorems}}
\Crefname{theorem}{\crefplainname{Theorem}}{\crefplainname{Theorems}}
\crefname{lemma}{\crefplainname{Lemma}}{\crefplainname{Lemmas}}
\Crefname{lemma}{\crefplainname{Lemma}}{\crefplainname{Lemmas}}
\crefname{proposition}{\crefplainname{Proposition}}{\crefplainname{Propositions}}
\Crefname{proposition}{\crefplainname{Proposition}}{\crefplainname{Propositions}}
\crefname{corollary}{\crefplainname{corollary}}{\crefplainname{corollaries}}
\Crefname{corollary}{\crefplainname{Corollary}}{\crefplainname{Corollaries}}
\crefname{assumption}{\crefplainname{Assumption}}{\crefplainname{Assumptions}}
\Crefname{assumption}{\crefplainname{Assumption}}{\crefplainname{Assumptions}}
\crefname{condition}{\crefplainname{Condition}}{\crefplainname{Conditions}}
\Crefname{condition}{\crefplainname{Condition}}{\crefplainname{Conditions}}
\crefname{section}{\crefplainname{Section}}{\crefplainname{Sections}}
\Crefname{section}{\crefplainname{Section}}{\crefplainname{Sections}}
\crefname{subsection}{\crefplainname{Section}}{\crefplainname{Sections}}
\Crefname{subsection}{\crefplainname{Section}}{\crefplainname{Sections}}
\crefname{subsubsection}{\crefplainname{Section}}{\crefplainname{Sections}}
\Crefname{subsubsection}{\crefplainname{Section}}{\crefplainname{Sections}}
\crefname{figure}{\crefplainname{figure}}{\crefplainname{figures}}
\Crefname{figure}{\crefplainname{Figure}}{\crefplainname{Figures}}
\crefname{table}{\crefplainname{table}}{\crefplainname{tables}}
\Crefname{table}{\crefplainname{Table}}{\crefplainname{Tables}}
\crefname{equation}{\crefplainname{equation}}{\crefplainname{equations}}
\Crefname{equation}{\crefplainname{Equation}}{\crefplainname{Equations}}

\makeatletter
\newcommand{\switchtosupplementnumbering}{%
  \clearpage
  \setcounter{section}{0}%
  \setcounter{subsection}{0}%
  \setcounter{subsubsection}{0}%
  \setcounter{equation}{0}%
  \setcounter{theorem}{0}%
  \setcounter{assumption}{0}%
  \setcounter{figure}{0}%
  \setcounter{table}{0}%
  \renewcommand{\thefigure}{S\arabic{figure}}%
  \renewcommand{\thetable}{S\arabic{table}}%
  \renewcommand{\thesection}{\Alph{section}}%
  \renewcommand{\thesubsection}{\thesection.\arabic{subsection}}%
  \renewcommand{\thesubsubsection}{\thesubsection.\arabic{subsubsection}}%
  \renewcommand{\theequation}{\thesection.\arabic{equation}}%
  \renewcommand{\thetheorem}{\thesection.\arabic{theorem}}%
  \renewcommand{\theassumption}{\thesection.\arabic{assumption}}%
  \numberwithin{equation}{section}%
  \numberwithin{assumption}{section}%
}
\makeatother
\begin{document}
\setlength{\bibsep}{0pt}
\begin{bibunit}[apalike]
\title{Inference for Fr\'echet Mean Sets}
\author{Ryunosuke Miyazaki\thanks{Email:
\href{mailto:miyazaki11051105@gmail.com}{miyazaki11051105@gmail.com}.}}
\date{Department of Social Data Science, Hitotsubashi University}
\maketitle
\begin{abstract}
We construct confidence sets for the entire Fr\'echet mean set by
inverting studentized pairwise risk inequalities. Gaussian multiplier
calibration localized around empirical near-minimizers gives asymptotic
coverage at least the nominal level, uniformly over distribution classes
satisfying empirical-process conditions. The construction requires
neither uniqueness nor a manifold structure. We characterize the limiting
binding directions and give conditions for exact asymptotic coverage.
On the intrinsic circle, exactness holds for every fixed non-point-mass
law and uniformly over a family in which two means merge at a quartic
minimum. With ideal multiplier quantiles, the construction reproduces
the Euclidean Wald set under both nonsingular and singular covariance.
Under regular local conditions, it is first-order equivalent to intrinsic
Wald inference.
\end{abstract}
\paperkeywords{Fr\'echet mean; metric space; multiplier bootstrap; random object; set-valued inference; Wald inference.}
\section{Introduction}
\label{sec:introduction}
Fr\'echet means minimize expected squared distance \citep{Frechet1948},
and the set of minimizers need not be a singleton. Existing work provides
local inference for regular unique means
\citep{BhattacharyaPatrangenaru2005,BhattacharyaLin2017}, empirical
likelihood on manifolds \citep{KurisuOtsu2026} and open books
\citep{BharathEtAl2024}, and Hausdorff-consistent
estimation of mean sets \citep{BlanchardJaffe2025}. We construct
confidence sets for the entire mean set, without requiring uniqueness
or imposing a specific geometric structure such as a manifold or an open book.

A candidate is a mean precisely when no competitor has smaller risk.
We invert these inequalities after studentization. The approach builds
on inference for minimizer sets
\citep{ChernozhukovHongTamer2007,RomanoShaikh2010,Bugni2010} and
generalized moment selection \citep{AndrewsSoares2010}. The sample-split
loss inversion of \citet{TakatsuKuchibhotla2025} controls coverage for
each minimizer separately; we calibrate a common multiplier critical
value for simultaneous inclusion of the entire mean set.

Even when the Fr\'echet mean is unique, comparisons with nearby
candidates can determine the asymptotic critical value.
As a candidate approaches the mean, the loss-difference variance
vanishes, but the centred and standardized loss difference can converge
to a nonzero function.
We characterize the limiting Gaussian supremum through $L^2(P)$ limits
of normalized contrasts between population means and competitors whose
risk gaps vanish relative to the contrast standard deviations.
The resulting class includes positive-variance comparisons between
distinct means and score directions at a regular unique mean.
Both global and localized multiplier calibration give uniform
asymptotic coverage on any distribution class satisfying our uniform
empirical-process conditions.
Under an additional condition for uniform recovery of the binding
directions, localized calibration consistently estimates the oracle
cutoff and attains the nominal level uniformly over that class.

The construction also retains the first-order precision of Wald
inference in regular models. With ideal multiplier quantiles, both
calibrations reproduce the Euclidean Wald set exactly, including under
structural rank deficiency. Under regular local conditions, the
localized set differs from the corresponding Wald set by
$o_p(n^{-1/2})$ in Hausdorff distance. Thus, extending inference to
possibly nonunique means does not entail first-order enlargement at
regular distributions. Under uniform empirical-process, risk-separation and margin conditions,
both confidence sets contract towards the mean set at a rate determined
by local risk growth, without using the margin exponent. Quadratic and
quartic margins yield confidence-radius bounds of order $n^{-1/2}$ and
$n^{-1/6}$, respectively.
On the intrinsic circle, localized coverage is asymptotically exact for
every fixed non-point-mass law, without density or curvature assumptions.
Coverage is also uniformly exact in a family where two means merge at a
quartic minimum.

\section{Studentized risk inversion and oracle calibration}
\label{sec:method}

\subsection{Confidence sets from pairwise risk inequalities}

Let $(\cQ,d)$ be a Polish metric space and let
$Y_1,\ldots,Y_n\stackrel{\mathrm{i.i.d}}{\sim}P$. Assume that
$\int d^2(y,o)\,dP(y)<\infty$ for some $o\in\cQ$.
Write $F_P(q)\coloneqq\int d^2(y,q)\,dP(y)$ and assume that
$M_P\coloneqq\argmin_{q\in\cQ}F_P(q)$ is nonempty.
We state uniform results for arbitrary nonempty classes $\mathcal P$
of such laws satisfying the assumptions of the corresponding result.
Each class is fixed as $n\to\infty$; all uniform bounds and limits are
taken over that class. The confidence sets do not depend on $\mathcal P$.
We use $\pr_P$ and $\E_P$ for sampling probability and expectation.
The $L^2(P)$ distance is
$\rho_P(f,g)\coloneqq\{\int(f-g)^2\,dP\}^{1/2}$; an overline denotes
closure in $L^2(P)$. We set $x_+\coloneqq\max(x,0)$ and
$\inf\varnothing=+\infty$. Empty suprema are zero for nonnegative
quantities and $-\infty$ otherwise.
For a real-valued random variable $Z$ and $0<p<1$, its $p$-quantile is
$q_p(Z)\coloneqq\inf\{t\in\mathbb R:\pr(Z\le t)\ge p\}$.

The empirical Fr\'echet function is
$F_n(q)\coloneqq n^{-1}\sum_{i=1}^n d^2(Y_i,q)$. The minimum population risk is
$F_P^\star\coloneqq\inf_{q\in\cQ}F_P(q)$, so that
$M_P=\{q\in\cQ:F_P(q)=F_P^\star\}$. Equivalently,
\[
 q\in M_P\quad\Longleftrightarrow\quad
 F_P(q)\le F_P(r)\quad\text{for every }r\in\cQ.
\]
For each pair $q,r$, define
\[
\begin{split}
 \sigma_P^2(q,r)&\coloneqq\int\bigl[d^2(y,q)-d^2(y,r)
       -\{F_P(q)-F_P(r)\}\bigr]^2\,dP(y),\\
 \wh\sigma_n^2(q,r)&\coloneqq\frac1n\sum_{i=1}^n
 \bigl[d^2(Y_i,q)-d^2(Y_i,r)
       -\{F_n(q)-F_n(r)\}\bigr]^2.
\end{split}
\]
Inverting the empirical risk inequalities gives a family of sets
indexed by a common nonnegative critical value:
\begin{equation}
 \mathcal C_n(c)\coloneqq\{q\in\cQ:n^{1/2}\{F_n(q)-F_n(r)\}
 \le c\wh\sigma_n(q,r)\ \text{for every }r\in\cQ\},\qquad c\ge0.
 \label{eq:inversion}
\end{equation}
The calibration problem is to choose a common $c$ for simultaneous
coverage of $M_P$. The procedures in \cref{sec:calibration} retain all
comparisons in \eqref{eq:inversion} and differ only in this cutoff.

\subsection{Oracle calibration and binding directions}

For $\wh\sigma_n(q,r)>0$, set
\[
 T_n(q,r)\coloneqq\frac{n^{1/2}\{F_n(q)-F_n(r)\}}{\wh\sigma_n(q,r)},
 \qquad
 T_n(M_P)\coloneqq\sup_{\substack{m\in M_P,\ r\in\cQ\\
                         \wh\sigma_n(m,r)>0}}T_n(m,r).
\]
To identify the pointwise asymptotic oracle cutoff, we determine which
comparisons contribute to the limit of $T_n(M_P)$.

For $\sigma_P(q,r)>0$, define the normalized contrast
\[
 h_{q,r}(y)\coloneqq\frac{d^2(y,q)-d^2(y,r)-\{F_P(q)-F_P(r)\}}
                       {\sigma_P(q,r)}.
\]
The class $\mathcal H_P\coloneqq
\{h_{q,r}:q,r\in\cQ,\ \sigma_P(q,r)>0\}$ is symmetric because
$h_{r,q}=-h_{q,r}$; each element has $P$-mean zero and $L^2(P)$ norm one.
For $h\in\mathcal H_P$, let
$\mathbb G_nh\coloneqq n^{-1/2}\sum_{i=1}^n h(Y_i)$.

For $m\in M_P$ and $\sigma_P(m,r)>0$, define the standardized
slackness by
\[
 \lambda_P(m,r)\coloneqq\frac{F_P(r)-F_P(m)}{\sigma_P(m,r)}\ge0.
\]
Whenever both population and empirical variances are positive,

\begin{equation}
 T_n(m,r)=\frac{\sigma_P(m,r)}{\wh\sigma_n(m,r)}
 \{\mathbb G_nh_{m,r}-n^{1/2}\lambda_P(m,r)\}.
 \label{eq:studentized-mean-decomposition}
\end{equation}

Under the regularity condition below, comparisons with standardized
slack bounded away from zero have diverging negative drift. The
limiting supremum is therefore governed by vanishing-slack directions.
Define the binding class by
\[
 \mathcal H_{P,0}\coloneqq\bigcap_{\eta>0}
 \overline{\{h_{m,r}:m\in M_P,\ r\in\cQ,\
 \sigma_P(m,r)>0,\ \lambda_P(m,r)\le\eta\}}^{\,L^2(P)}.
\]
It includes positive-variance comparisons between distinct means and
limiting local directions, which may arise even when the mean is unique.

\begin{assumption}[Uniform Donsker regularity on $\mathcal P$]
\label{ass:main-process}
For the specified nonempty distribution class $\mathcal P$,
$\mathcal H_P$ is nonempty for every $P\in\mathcal P$.
The family $\{\mathcal H_P:P\in\mathcal P\}$ is Donsker and
pre-Gaussian uniformly over $\mathcal P$.
There are measurable envelopes $H_P$ satisfying
\(
 \lim_{L\to\infty}\sup_{P\in\mathcal P}
 \int H_P^2\mathds 1\{H_P>L\}\,dP=0.
\)
\end{assumption}
For fixed-law results, \cref{ass:main-process} is applied with
$\mathcal P=\{P\}$. A nonempty $P$-Donsker class $\mathcal H_P$
with a measurable square-integrable envelope suffices in this case.
The uniform framework \citep{SheehyWellner1992,VaartWellner2023}
is specified in \suppcref{sec:supp-studentization}, and
\suppcref{sec:supp-entropy} provides sufficient conditions, including finite
centred-loss rank and finite-dimensional H\"older classes.

Whenever \cref{ass:main-process} holds for $\mathcal P=\{P\}$,
the centred Gaussian limit $\mathbb G_P$ has covariance
$\Cov(\mathbb G_P(h),\mathbb G_P(g))=\int h(y)g(y)\,dP(y)$
and is extended continuously to
$\overline{\mathcal H_P}^{\,L^2(P)}$, which contains $\mathcal H_{P,0}$.

\begin{theorem}[Binding Gaussian limit]
\label{thm:binding-main}
\samepage
Fix a law $P$ and suppose that \cref{ass:main-process} holds
for $\mathcal P=\{P\}$. If $\mathcal H_{P,0}\ne\varnothing$, then
\begin{equation*}
 T_n(M_P)\rightsquigarrow T_P^0
 \coloneqq\sup_{h\in\mathcal H_{P,0}}\mathbb G_P(h).
\end{equation*}
If $\mathcal H_{P,0}=\varnothing$, then $T_n(M_P)\to_p-\infty$.
\end{theorem}

For $0<\alpha<1/2$ and a nonempty binding class, define
$c_{P,\alpha}^0\coloneqq q_{1-\alpha}(T_P^0)$, the pointwise
asymptotic oracle cutoff identified by \cref{thm:binding-main}.
Set $c_{P,\alpha}^0\coloneqq0$ when the binding class is empty.

\section{Uniform coverage and localized calibration}
\label{sec:calibration}

\subsection{Global calibration}

Let $\xi_1,\ldots,\xi_n$ be independent standard normals,
independent of the data. For $\wh\sigma_n(q,r)>0$, define
\[
 T_n^\#(q,r)\coloneqq\frac{n^{-1/2}\sum_{i=1}^n\xi_i
 [d^2(Y_i,q)-d^2(Y_i,r)-\{F_n(q)-F_n(r)\}]}
 {\wh\sigma_n(q,r)}.
\]
Global calibration uses all positive-variance comparisons. Set
$S_n^{\#,G}\coloneqq(\sup_{\wh\sigma_n(q,r)>0}T_n^\#(q,r))_+$,
let $c_n^G$ be its conditional $(1-\alpha)$-quantile, and define
$C_n^G\coloneqq\mathcal C_n(c_n^G)$. When $\mathcal H_P\ne\varnothing$, its limiting
target is
$T_P^G\coloneqq\sup_{h\in\mathcal H_P}\mathbb G_P(h)
=\sup_{h\in\mathcal H_P}|\mathbb G_P(h)|$.
\begin{theorem}[Uniform global calibration]
\label{thm:calibration}
Fix $0<\alpha<1/2$. Let $\mathcal P$ be any nonempty class
satisfying \cref{ass:main-process}. Then
\begin{equation*}
 \liminf_{n\to\infty}\inf_{P\in\mathcal P}
 \pr_P(M_P\subseteq C_n^G)\ge1-\alpha.
\end{equation*}
Moreover, for every $\epsilon>0$,
\[
 \sup_{P\in\mathcal P}
 \pr_P\bigl(|c_n^G-c_{P,\alpha}^G|>\epsilon\bigr)\to0,
 \qquad c_{P,\alpha}^G\coloneqq q_{1-\alpha}(T_P^G)>0.
\]
\end{theorem}
Since $\mathcal H_{P,0}\subseteq\overline{\mathcal H_P}^{\,L^2(P)}$,
we have $c_{P,\alpha}^0\le c_{P,\alpha}^G$ whenever
$\mathcal H_{P,0}\ne\varnothing$, with strict inequality under strict
inclusion (\suppcref{lem:gaussian-quantile}).

\subsection{Localized calibration and exactness}

To target binding directions, select comparisons anchored at
empirical near-minimizers and having small empirical slack.
The population near-minimizer set is
$M_P(\delta)\coloneqq\{q:F_P(q)-F_P^\star\le\delta\}$ for
$\delta\ge0$, and the minimum empirical risk is
$F_n^\star\coloneqq\inf_{s\in\cQ}F_n(s)$.
Let $b_n,\tau_n>0$ be deterministic sequences satisfying
\begin{equation}
 b_n\to\infty,\qquad
 b_nn^{-1/2}\to0,\qquad
 \tau_n\to\infty,\qquad
 \tau_n n^{-1/2}\to0.
 \label{eq:adaptive-tuning-main}
\end{equation}
Define the empirical relative near-minimizer set by
\[
 \wh M_n\coloneqq\left\{q\in\cQ:
 F_n(q)\le\left(1+n^{-1/2}b_n\right)F_n^\star\right\}.
\]
The relative anchor is nonempty for every sample.
Under \cref{ass:main-process}, it contains $M_P$ with probability
tending to one uniformly over the specified class $\mathcal P$,
without requiring $\sup_{P\in\mathcal P}F_P^\star<\infty$
(\suppcref{lem:bounded-ball-localization}).

Define the localized selector directly by
\[
 \widehat{\mathcal B}_n^{\mathrm L}\coloneqq\left\{(q,r):
 q\in\wh M_n,\ \widehat\sigma_n(q,r)>0,\ T_n(q,r)\ge-\tau_n\right\}.
\]
Set
\[
 S_n^{\#,\mathrm L}\coloneqq\left(\sup_{(q,r)\in\wh{\mathcal B}_n^{\mathrm L}}
 T_n^\#(q,r)\right)_+,\quad
 c_n^{\mathrm L}\coloneqq q_{1-\alpha}(S_n^{\#,\mathrm L}\mid Y_1,\ldots,Y_n),
 \quad C_n^{\mathrm L}\coloneqq\mathcal C_n(c_n^{\mathrm L}).
\]
An empty selection has $S_n^{\#,\mathrm L}=0$.
Monotonicity of suprema, conditional quantiles and inversion gives
$c_n^{\mathrm L}\le c_n^G$ and $C_n^{\mathrm L}\subseteq C_n^G$.
Define the low-risk secant class
\[
 \mathcal S_P(\delta)\coloneqq
 \overline{\{h_{q,r}:q,r\in M_P(\delta),\ \sigma_P(q,r)>0\}}^{\,L^2(P)}.
\]
For $\mathcal H_{P,0}\ne\varnothing$ and $\delta>0$, define the directed
$L^2(P)$ excess of $\mathcal S_P(\delta)$ over the binding class by
\[
 \omega_P(\delta)\coloneqq
 \sup_{g\in\mathcal S_P(\delta)}
 \inf_{h\in\mathcal H_{P,0}}\rho_P(g,h).
\]
\begin{theorem}[Uniform localized calibration]
\label{thm:zero-cutoff-main}
Fix $0<\alpha<1/2$ and assume \eqref{eq:adaptive-tuning-main}.
Let $\mathcal P$ be any nonempty class satisfying
\cref{ass:main-process}. Then
\begin{equation*}
 \liminf_n\inf_{P\in\mathcal P}
 \pr_P(M_P\subseteq C_n^{\mathrm L})\ge1-\alpha.
\end{equation*}
Let $\mathcal P_s\subseteq\mathcal P$ be any nonempty subclass
such that $\mathcal H_{P,0}\ne\varnothing$ for every
$P\in\mathcal P_s$ and
\begin{equation}
 \lim_{u\downarrow0}\sup_{P\in\mathcal P_s}
 \omega_P(uF_P^\star)=0.
 \label{eq:uniform-secants}
\end{equation}
Then, for every $\epsilon>0$,
\begin{equation*}
 \sup_{P\in\mathcal P_s}
 \pr_P(|c_n^{\mathrm L}-c_{P,\alpha}^0|>\epsilon)\to0,\qquad
 \sup_{P\in\mathcal P_s}
 \left|\pr_P(M_P\subseteq C_n^{\mathrm L})-(1-\alpha)\right|\to0.
\end{equation*}
\end{theorem}
Condition \eqref{eq:uniform-secants} is required only on the subclass
for which oracle-cutoff consistency and exactness are asserted.
The uniform coverage bound does not require this condition.
Sufficient conditions are given in \suppcref{sec:supp-directions}.

\begin{corollary}[Intrinsic circle]
\label{cor:circle-exactness-main}
Let $\cQ=S^1$ have its intrinsic geodesic distance. Fix
$0<\alpha<1/2$ and assume \eqref{eq:adaptive-tuning-main}.
For every fixed non-point-mass Borel law $P$ on $S^1$,
\[
 c_n^{\mathrm L}\to_p c_{P,\alpha}^0,\qquad
 \pr_P(M_P\subseteq C_n^{\mathrm L})\longrightarrow1-\alpha.
\]
Both convergences hold uniformly over the family
$\{P_\lambda:-1/8\le\lambda\le1\}$ in \eqref{eq:circle-family}.
\end{corollary}
The uniform assertion includes the merger of two means at $\lambda=0$,
where the risk has a quartic minimum.

\begin{remark}[Pointwise refinements and uniform validity]
\label{subsec:pointwise-refinements}
Additional restrictions on the calibration directions can reduce
conservativeness and recover the binding critical value at a fixed
distribution under suitable separation and retention conditions.
\suppCref{sec:pointwise-refinements} gives the scale-restricted
construction and its pointwise guarantees; a tripod calculation
there separately identifies a constraint on uniform oracle
adaptation for the common scalar inversion \eqref{eq:inversion}.
\end{remark}

\section{Wald recovery and contraction}
\label{sec:geometry}
Coverage and contraction rates do not determine the first-order size
of a confidence set. We therefore compare the proposed construction
with Wald inference. The first result gives exact equality in
Euclidean spaces, while the second establishes local equivalence
under regularity. The final result gives confidence-radius bounds
under more general risk margins.

First let $\cQ\coloneqq\mathbb R^d$, with fixed $d$,
$\mu\coloneqq\int y\,dP(y)$ and
$\Sigma\coloneqq\int(y-\mu)(y-\mu)^{\mathrm T}\,dP(y)\succeq0$, with
$r\coloneqq\operatorname{rank}(\Sigma)\ge1$.
The identity
$d^2(y,q)-d^2(y,s)-\{F_P(q)-F_P(s)\}
=2(s-q)^{\mathrm T}(y-\mu)$
shows that every global normalized direction is a normalized linear
projection on $\operatorname{Range}(\Sigma)$. Each is also the limit of a
comparison between $\mu$ and $\mu+tv$, for $v\in\operatorname{Range}(\Sigma)
\setminus\{0\}$, as $t\downarrow0$; hence
$\mathcal H_P=\mathcal H_{P,0}$.

\begin{proposition}[Finite-sample Euclidean reduction]
\label{prop:euclidean}
Fix $0<\alpha<1/2$. Put $\bar Y_n\coloneqq n^{-1}\sum_iY_i$ and
$\wh\Sigma_n\coloneqq n^{-1}\sum_i(Y_i-\bar Y_n)(Y_i-\bar Y_n)^{\mathrm T}$.
Let $\wh r_n\coloneqq\operatorname{rank}(\wh\Sigma_n)$.
Write $A^+$ for the Moore--Penrose inverse of a matrix $A$.
\begin{enumerate}
\item[(i)] On a sample with $\wh r_n>0$, the supremum over
positive-variance comparisons satisfies
\begin{equation}
 \sup_{s:\wh\sigma_n(q,s)>0}T_n(q,s)
 =\begin{cases}
 \left\{n(q-\bar Y_n)^{\mathrm T}\wh\Sigma_n^+(q-\bar Y_n)\right\}^{1/2},
   &q-\bar Y_n\in\operatorname{Range}(\wh\Sigma_n),\\
 +\infty,&q-\bar Y_n\notin\operatorname{Range}(\wh\Sigma_n).
 \end{cases}
 \label{eq:euclidean-statistic}
\end{equation}

With all comparators and ideal multiplier quantiles,
both calibrations satisfy
\begin{equation}
C_n^G=C_n^{\mathrm L}=
 \left\{q\in\bar Y_n+\operatorname{Range}(\wh\Sigma_n):
 n(q-\bar Y_n)^{\mathrm T}\wh\Sigma_n^+(q-\bar Y_n)
 \leq\chi^2_{\wh r_n,1-\alpha}\right\}.
 \label{eq:euclidean}
\end{equation}
If $\wh r_n=0$, both cutoffs are zero and both sets equal
$\{\bar Y_n\}$.
\item[(ii)]
Let $\mathcal P$ be any nonempty class of laws on $\mathbb R^d$
with finite second moments, means $\mu_P$, covariances $\Sigma_P$
and ranks $r_P\coloneqq\operatorname{rank}(\Sigma_P)\ge1$,
satisfying \cref{ass:main-process}.
Both Gaussian suprema in \cref{thm:binding-main,thm:calibration}
have law $\chi_{r_P}$, and
\[
 \inf_{P\in\mathcal P}\pr_P(\wh r_n=r_P)\to1,\qquad
 \sup_{P\in\mathcal P}\max_{J\in\{G,\mathrm L\}}
 \left|\pr_P(\mu_P\in C_n^J)-(1-\alpha)\right|\to0.
\]

\end{enumerate}
\end{proposition}

On $\{\wh r_n=r\}$, \eqref{eq:euclidean} is the ordinary Wald ellipsoid
if $r=d$ and its lower-dimensional generalized Wald analogue on
$\bar Y_n+\operatorname{Range}(\wh\Sigma_n)$ if $r<d$. The pairwise construction yields both the
support constraint and rank-adapted calibration without taking the rank
or a generalized inverse as an input. This reduction concerns structural
rank deficiency.

For local Wald recovery in a general metric space, use the chart $\phi$ and empirical mean
coordinate $\widehat u_n$ of \suppcref{sec:supp-wald-conditions}.
The local risk curvature $H$ and score covariance $\Sigma$ are positive
definite, with $V\coloneqq H^{-1}\Sigma H^{-1}$. For a consistent
covariance estimate $\widehat V_n$, define
\begin{equation*}
 \mathcal E_n(a)\coloneqq\left\{\phi(u):u\in B(0,r_0),\quad
 n(u-\widehat u_n)^{\mathrm T}\widehat V_n^{-1}(u-\widehat u_n)
 \leq a^2\right\},\qquad a\geq0.
\end{equation*}

Here $d_{H,d}$ denotes Hausdorff distance induced by $d$.

\begin{proposition}[Local Wald recovery]
\label{prop:wald-recovery}
Fix $P$ and $0<\alpha<1/2$. Assume the regular local conditions of
\suppcref{sec:supp-wald-conditions} and
\eqref{eq:adaptive-tuning-main}. For
$c_\alpha\coloneqq(\chi^2_{d,1-\alpha})^{1/2}$,
\begin{equation*}
 c_n^{\mathrm L}=c_\alpha+o_p(1),\qquad
 n^{1/2}d_{H,d}(C_n^{\mathrm L},\mathcal E_n(c_\alpha))=o_p(1),
\end{equation*}
and $\pr_P(M_P\subseteq C_n^{\mathrm L})\to1-\alpha$.
\end{proposition}

The Hausdorff error is smaller than the $n^{-1/2}$ scale of the
regular confidence set. In local coordinates, the construction
therefore recovers the first-order centre, shape and directional
widths of the Wald set. Within the inversion family \eqref{eq:inversion},
any nonnegative common cutoff converging to a constant smaller than
$c_\alpha$ gives asymptotic undercoverage at this fixed distribution
(\suppcref{sec:supp-wald-proof}).

If \cref{ass:main-process} also holds for $\mathcal P=\{P\}$,
global calibration has the same centre and covariance shape but
radius $c_{P,\alpha}^G\ge c_\alpha$, which may be larger.

For $A\subseteq\cQ$ and nonempty
$B\subseteq\cQ$, let
$d(q,B)\coloneqq\inf_{b\in B}d(q,b)$ and
$\operatorname{rad}(A;B)\coloneqq\sup_{q\in A}d(q,B)$, with radius zero
for empty $A$.
\begin{proposition}[Uniform adaptive confidence radius]
\label{prop:adaptive-radius}
Let $\mathcal P$ be any nonempty class satisfying
\cref{ass:main-process}. Let $c_n\ge0$ be measurable with
$c_n=O_{\mathcal P}(1)$. Suppose that:
\begin{enumerate}
\item[(i)] the minimum risks satisfy
$\overline F\coloneqq\sup_{P\in\mathcal P}F_P^\star<\infty$;

\item[(ii)] for every $\epsilon>0$,

\[
 \inf_{P\in\mathcal P}\inf_{q:\ d(q,M_P)\ge\epsilon}
 \{F_P(q)-F_P^\star\}>0;
\]

\item[(iii)] there are constants
$a,r_0>0$ and $\beta>1$, which may depend on $\mathcal P$
but not on $P\in\mathcal P$ or $n$, such that, for every
$P\in\mathcal P$ and every $q\in\cQ$ satisfying
$0<d(q,M_P)\le r_0$,
\begin{equation}
 F_P(q)-F_P^\star\ge a\,d(q,M_P)^\beta.
 \label{eq:margin-scale-main}
\end{equation}
\end{enumerate}
Then
\begin{equation}
 \operatorname{rad}(\mathcal C_n(c_n);M_P)
 =O_{\mathcal P}\bigl(n^{-1/\{2(\beta-1)\}}\bigr).
 \label{eq:adaptive-radius-main}
\end{equation}
The conclusion applies to both global and localized calibration.
\end{proposition}
The margin condition \citep{Schoetz2019} includes quadratic barycentric
variance inequalities \citep{AhidarCoutrixEtAl2020}; higher-order
flatness includes smeary means \citep{EltznerHuckemann2019}.
Squared-distance loss has linear contrast scale
(\suppcref{lem:automatic-linear-scale}), so no separate scale condition
is needed. The set does not use $\beta$, and its radius adapts to the
risk margin. Quadratic and quartic margins give rates $n^{-1/2}$ and
$n^{-1/6}$, respectively.

\section{Numerical experiments}
\label{sec:numerics}
\begin{figure}
\centering
\paperfigurewidth
\begin{tikzpicture}
\begin{scope}[yshift=0cm]
\begin{groupplot}[paperaxis,width=3.9cm,group style={group size=3 by 1,horizontal sep=.6cm}]
\nextgroupplot[title={Euclidean: $\delta=0$},ylabel={Coverage (\%)},]
\addplot[black!55,densely dashed,line width=.45pt] coordinates {(40,95)(250,95)};
\addplot[paperg,error bars/y dir=both,error bars/y explicit]
table[x=x,y=y,y error minus=lo,y error plus=up,row sep=\\] {
x y lo up\\
48.886228 92.400000 1.820432 1.565634\\
97.772456 95.500000 1.475283 1.199037\\
195.544912 94.600000 1.587401 1.317665\\
};
\addplot[paperwald,error bars/y dir=both,error bars/y explicit]
table[x=x,y=y,y error minus=lo,y error plus=up,row sep=\\] {
x y lo up\\
51.139147 92.300000 1.829938 1.575797\\
102.278294 95.500000 1.475283 1.199037\\
204.556588 94.500000 1.599152 1.330120\\
};
\nextgroupplot[title={Euclidean: $\delta=1/2$},yticklabels={},]
\addplot[black!55,densely dashed,line width=.45pt] coordinates {(40,95)(250,95)};
\addplot[paperg,error bars/y dir=both,error bars/y explicit]
table[x=x,y=y,y error minus=lo,y error plus=up,row sep=\\] {
x y lo up\\
48.886228 92.400000 1.820432 1.565634\\
97.772456 94.100000 1.644907 1.378665\\
195.544912 94.500000 1.599152 1.330120\\
};
\addplot[paperwald,error bars/y dir=both,error bars/y explicit]
table[x=x,y=y,y error minus=lo,y error plus=up,row sep=\\] {
x y lo up\\
51.139147 92.500000 1.810845 1.555387\\
102.278294 94.000000 1.656049 1.390496\\
204.556588 94.200000 1.633650 1.366714\\
};
\nextgroupplot[title={Euclidean: $\delta=1$},yticklabels={},]
\addplot[black!55,densely dashed,line width=.45pt] coordinates {(40,95)(250,95)};
\addplot[paperg,error bars/y dir=both,error bars/y explicit]
table[x=x,y=y,y error minus=lo,y error plus=up,row sep=\\] {
x y lo up\\
48.886228 92.400000 1.820432 1.565634\\
97.772456 94.800000 1.563501 1.292347\\
195.544912 92.900000 1.771644 1.513535\\
};
\addplot[paperwald,error bars/y dir=both,error bars/y explicit]
table[x=x,y=y,y error minus=lo,y error plus=up,row sep=\\] {
x y lo up\\
51.139147 92.600000 1.801174 1.545056\\
102.278294 95.100000 1.526605 1.253294\\
204.556588 93.100000 1.751511 1.492066\\
};
\end{groupplot}
\end{scope}
\begin{scope}[yshift=-3.5cm]
\begin{groupplot}[paperaxis,width=2.75cm,group style={group size=4 by 1,horizontal sep=.6cm}]
\nextgroupplot[title={Ranking: Euclidean},ylabel={Coverage (\%)},]
\addplot[black!55,densely dashed,line width=.45pt] coordinates {(40,95)(250,95)};
\addplot[paperg,error bars/y dir=both,error bars/y explicit]
table[x=x,y=y,y error minus=lo,y error plus=up,row sep=\\] {
x y lo up\\
48.886228 96.800000 1.287590 1.001115\\
97.772456 98.300000 1.007950 0.706649\\
195.544912 97.000000 1.255140 0.966951\\
};
\addplot[paperl,error bars/y dir=both,error bars/y explicit]
table[x=x,y=y,y error minus=lo,y error plus=up,row sep=\\] {
x y lo up\\
51.139147 93.900000 1.667077 1.402212\\
102.278294 94.900000 1.551346 1.279476\\
204.556588 93.700000 1.688804 1.425307\\
};
\nextgroupplot[title={Ranking: Kendall},yticklabels={},]
\addplot[black!55,densely dashed,line width=.45pt] coordinates {(40,95)(250,95)};
\addplot[paperg,error bars/y dir=both,error bars/y explicit]
table[x=x,y=y,y error minus=lo,y error plus=up,row sep=\\] {
x y lo up\\
48.886228 99.500000 0.662947 0.337458\\
97.772456 99.600000 0.620966 0.290909\\
195.544912 99.600000 0.620966 0.290909\\
};
\addplot[paperl,error bars/y dir=both,error bars/y explicit]
table[x=x,y=y,y error minus=lo,y error plus=up,row sep=\\] {
x y lo up\\
51.139147 93.900000 1.667077 1.402212\\
102.278294 95.900000 1.421246 1.141983\\
204.556588 94.400000 1.610775 1.342445\\
};
\nextgroupplot[title={Ranking: footrule},yticklabels={},]
\addplot[black!55,densely dashed,line width=.45pt] coordinates {(40,95)(250,95)};
\addplot[paperg,error bars/y dir=both,error bars/y explicit]
table[x=x,y=y,y error minus=lo,y error plus=up,row sep=\\] {
x y lo up\\
48.886228 99.200000 0.770205 0.454002\\
97.772456 99.800000 0.520584 0.175770\\
195.544912 98.900000 0.859663 0.549641\\
};
\addplot[paperl,error bars/y dir=both,error bars/y explicit]
table[x=x,y=y,y error minus=lo,y error plus=up,row sep=\\] {
x y lo up\\
51.139147 94.200000 1.633650 1.366714\\
102.278294 95.200000 1.514012 1.239973\\
204.556588 93.300000 1.731004 1.470216\\
};
\nextgroupplot[title={Ranking: graph},yticklabels={},]
\addplot[black!55,densely dashed,line width=.45pt] coordinates {(40,95)(250,95)};
\addplot[paperg,error bars/y dir=both,error bars/y explicit]
table[x=x,y=y,y error minus=lo,y error plus=up,row sep=\\] {
x y lo up\\
48.886228 98.100000 1.051240 0.752296\\
97.772456 98.900000 0.859663 0.549641\\
195.544912 98.900000 0.859663 0.549641\\
};
\addplot[paperl,error bars/y dir=both,error bars/y explicit]
table[x=x,y=y,y error minus=lo,y error plus=up,row sep=\\] {
x y lo up\\
51.139147 94.800000 1.563501 1.292347\\
102.278294 96.400000 1.349291 1.066109\\
204.556588 95.900000 1.421246 1.141983\\
};
\end{groupplot}
\end{scope}
\begin{scope}[yshift=-7cm]
\begin{groupplot}[paperaxis,width=2.75cm,height=2cm,
 ymin=88,ymax=101,ytick={90,95,100},xtick=\empty,
 clip marker paths=true,
 axis x line=none,
 group style={group name=circleupper,group size=4 by 1,horizontal sep=.6cm}]
\nextgroupplot[title={Circle: $\lambda=-1/8$},ylabel={Coverage (\%)},]
\addplot[black!55,densely dashed,line width=.45pt] coordinates {(40,95)(250,95)};
\addplot[paperg,error bars/y dir=both,error bars/y explicit]
table[x=x,y=y,y error minus=lo,y error plus=up,row sep=\\] {
x y lo up\\
46.976137 98.200000 1.029904 0.729805\\
93.952275 99.300000 0.736919 0.418114\\
187.904550 98.900000 0.859663 0.549641\\
};
\addplot[paperl,error bars/y dir=both,error bars/y explicit]
table[x=x,y=y,y error minus=lo,y error plus=up,row sep=\\] {
x y lo up\\
48.971015 93.400000 1.720605 1.459143\\
97.942030 95.700000 1.448626 1.170883\\
195.884060 96.500000 1.334238 1.050247\\
};
\addplot[paperel,error bars/y dir=both,error bars/y explicit]
table[x=x,y=y,y error minus=lo,y error plus=up,row sep=\\] {
x y lo up\\
51.050606 89.600000 2.060016 1.823218\\
102.101213 90.500000 1.988484 1.745979\\
204.202425 91.400000 1.912027 1.663747\\
};

\nextgroupplot[title={Circle: $\lambda=0$},yticklabels={},]
\addplot[black!55,densely dashed,line width=.45pt] coordinates {(40,95)(250,95)};
\addplot[paperg,error bars/y dir=both,error bars/y explicit]
table[x=x,y=y,y error minus=lo,y error plus=up,row sep=\\] {
x y lo up\\
46.976137 99.300000 0.736919 0.418114\\
93.952275 99.800000 0.520584 0.175770\\
187.904550 99.700000 0.574202 0.238090\\
};
\addplot[paperl,error bars/y dir=both,error bars/y explicit]
table[x=x,y=y,y error minus=lo,y error plus=up,row sep=\\] {
x y lo up\\
48.971015 95.800000 1.435029 1.156529\\
97.942030 97.500000 1.168482 0.875746\\
195.884060 98.400000 0.985325 0.682768\\
};
\addplot[paperel,error bars/y dir=both,error bars/y explicit]
table[x=x,y=y,y error minus=lo,y error plus=up,row sep=\\] {
x y lo up\\
51.050606 94.400000 1.610775 1.342445\\
102.101213 94.700000 1.575518 1.305075\\
204.202425 94.900000 1.551346 1.279476\\
};

\nextgroupplot[title={Circle: $\lambda=1/2$},yticklabels={},]
\addplot[black!55,densely dashed,line width=.45pt] coordinates {(40,95)(250,95)};
\addplot[paperg,error bars/y dir=both,error bars/y explicit]
table[x=x,y=y,y error minus=lo,y error plus=up,row sep=\\] {
x y lo up\\
46.976137 99.400000 0.701342 0.379502\\
93.952275 99.800000 0.520584 0.175770\\
187.904550 99.900000 0.455892 0.097468\\
};
\addplot[paperl,error bars/y dir=both,error bars/y explicit]
table[x=x,y=y,y error minus=lo,y error plus=up,row sep=\\] {
x y lo up\\
48.971015 97.000000 1.255140 0.966951\\
97.942030 97.800000 1.111966 0.816259\\
195.884060 99.100000 0.801578 0.487660\\
};
\addplot[paperel,error bars/y dir=both,error bars/y explicit]
table[x=x,y=y,y error minus=lo,y error plus=up,row sep=\\] {
x y lo up\\
51.050606 94.200000 1.633650 1.366714\\
102.101213 94.500000 1.599152 1.330120\\
204.202425 94.000000 1.656049 1.390496\\
};
\addplot[paperwald,error bars/y dir=both,error bars/y explicit]
table[x=x,y=y,y error minus=lo,y error plus=up,row sep=\\] {
x y lo up\\
212.874036 86.500000 2.275592 2.058095\\
};
\nextgroupplot[title={Circle: $\lambda=1$},yticklabels={},]
\addplot[black!55,densely dashed,line width=.45pt] coordinates {(40,95)(250,95)};
\addplot[paperg,error bars/y dir=both,error bars/y explicit]
table[x=x,y=y,y error minus=lo,y error plus=up,row sep=\\] {
x y lo up\\
46.976137 99.400000 0.701342 0.379502\\
93.952275 99.600000 0.620966 0.290909\\
187.904550 99.800000 0.520584 0.175770\\
};
\addplot[paperl,error bars/y dir=both,error bars/y explicit]
table[x=x,y=y,y error minus=lo,y error plus=up,row sep=\\] {
x y lo up\\
48.971015 96.900000 1.271509 0.984183\\
97.942030 98.400000 0.985325 0.682768\\
195.884060 98.600000 0.937816 0.632550\\
};
\addplot[paperel,error bars/y dir=both,error bars/y explicit]
table[x=x,y=y,y error minus=lo,y error plus=up,row sep=\\] {
x y lo up\\
51.050606 93.500000 1.710107 1.447969\\
102.101213 94.200000 1.633650 1.366714\\
204.202425 93.800000 1.677995 1.413815\\
};
\addplot[paperwald,error bars/y dir=both,error bars/y explicit]
table[x=x,y=y,y error minus=lo,y error plus=up,row sep=\\] {
x y lo up\\
53.218509 93.200000 1.741305 1.481190\\
106.437018 94.300000 1.622273 1.354642\\
212.874036 94.100000 1.644907 1.378665\\
};
\end{groupplot}
\foreach \p in {1,2,3,4} {
 \draw[black!65,line width=.45pt]
 ([shift={(-2pt,-2pt)}]circleupper c\p r1.south west) -- ++(4pt,2pt);
 \draw[black!65,line width=.45pt]
 ([shift={(-2pt,-4pt)}]circleupper c\p r1.south west) -- ++(4pt,2pt);
}
\begin{scope}[yshift=-.85cm]
\begin{groupplot}[paperaxis,width=2.75cm,height=.7cm,
 ymin=0,ymax=88,ytick={0,50},clip marker paths=false,
 group style={group size=4 by 1,horizontal sep=.6cm}]
\nextgroupplot[xlabel={$n$},ylabel={},]

\addplot[paperel,error bars/y dir=both,error bars/y explicit]
table[x=x,y=y,y error minus=lo,y error plus=up,row sep=\\] {
x y lo up\\
51.050606 89.600000 2.060016 1.823218\\
};
\addplot[paperwald,error bars/y dir=both,error bars/y explicit]
table[x=x,y=y,y error minus=lo,y error plus=up,row sep=\\] {
x y lo up\\
53.218509 29.400000 2.809629 2.931472\\
106.437018 22.700000 2.562810 2.724600\\
212.874036 17.600000 2.312171 2.504651\\
};
\nextgroupplot[xlabel={$n$},yticklabels={},]

\addplot[paperwald,error bars/y dir=both,error bars/y explicit]
table[x=x,y=y,y error minus=lo,y error plus=up,row sep=\\] {
x y lo up\\
53.218509 66.400000 3.022538 2.925614\\
106.437018 59.000000 3.121165 3.068021\\
212.874036 52.800000 3.148496 3.131968\\
};
\nextgroupplot[xlabel={$n$},yticklabels={},]

\addplot[paperwald,error bars/y dir=both,error bars/y explicit]
table[x=x,y=y,y error minus=lo,y error plus=up,row sep=\\] {
x y lo up\\
53.218509 83.600000 2.443133 2.243372\\
106.437018 83.200000 2.464106 2.266775\\
212.874036 86.500000 2.275592 2.058095\\
};
\nextgroupplot[xlabel={$n$},yticklabels={},]

\end{groupplot}
\end{scope}
\end{scope}
\begin{scope}[yshift=3cm]
\node[anchor=west,font=\fontsize{8}{9}\selectfont] at (-.5,0) {\textcolor{paperGlobal}{$\bullet$}\ Global\quad\textcolor{paperLocal}{$\blacksquare$}\ Local\quad\textcolor{paperEL}{$\diamond$}\ EL\quad$\triangle$\ Wald};
\end{scope}
\end{tikzpicture}
\caption{Whole-mean-set coverage (\%) with individual 95\% exact binomial
Monte Carlo intervals (1,000 datasets per setting). Dashed lines mark 95\%,
also the exact Gaussian reference in the Euclidean panels.}
\label{fig:coverage}
\end{figure}
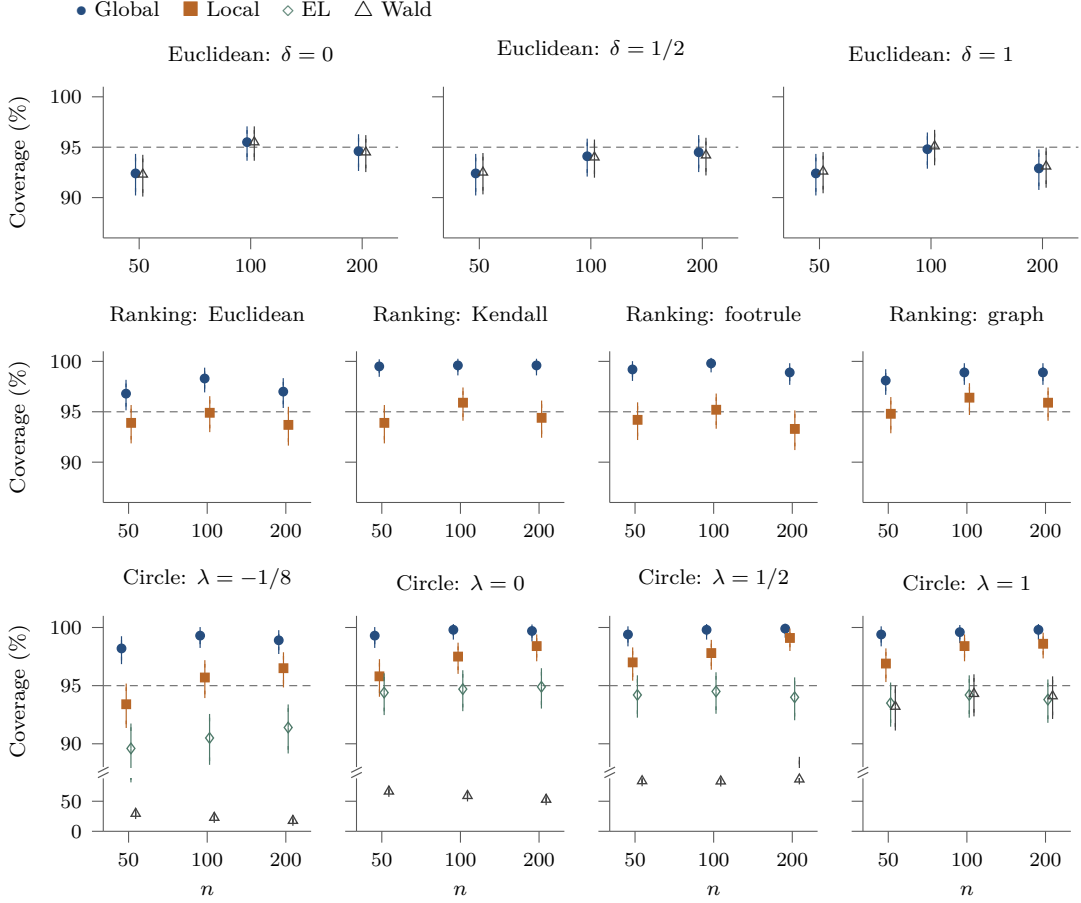

For each sample size $n\in\{50,100,200\}$, we generate 1,000 independent
datasets per setting and use nominal coverage 0.95, $B=1999$ Gaussian
multiplier draws, $b_n=(\log\log n)^{1/2}$ and $\tau_n=n^{1/5}$. Global and Local use $C_n^G$ and
$C_n^{\mathrm L}$, respectively, with shared observations and draws.
\Cref{fig:coverage} reports whole-mean-set coverage;
\suppcref{sec:supp-numerics} gives computational details.

Euclidean observations follow $Y\sim N_3(0,I_3-uu^{\mathrm T}+\delta uu^{\mathrm T})$,
where $u=(1,1,1)^{\mathrm T}/3^{1/2}$ and $\delta=0,1/2,1$.
At $\delta=0$, computations use orthonormal coordinates
on $u^\perp$, shared by all methods.
Global and Local coincide samplewise. Coverage is 92.4--95.5\%, within
0.3 percentage points of asymptotic $\chi^2$ Wald calibration.

Rankings follow $P=\tfrac12P_0+\tfrac12U_{24}$, where $P_0$ is uniform on
$\{1234,2314,3124\}$ and $U_{24}$ on all complete rankings.
On the 24-point candidate space, use rank-vector Euclidean, Kendall
and footrule distances \citep{Marden1995,DiaconisGraham1977}.
The squared losses give six, three and three means, with different
Kendall and footrule mean sets.
Global coverage is conservative; Local coverage is 93.3--95.9\%.
For footrule at $n=200$, even the Monte Carlo upper endpoint is below 95\%.
The continuous adjacent-swap graph has six edge-interior means
(\suppcref{sec:supp-population}) and coverage 98.1--98.9\% for Global,
94.8--96.4\% for Local.

For circular observations, use $P_\lambda$
with arc-length density
\begin{equation}
 f_\lambda(x)=\frac{1+(4-\lambda)\cos x/6
                    +2(1-\lambda)\cos(2x)/3}{2\pi},
 \qquad x\in[-\pi,\pi),\quad -1/8\le\lambda\le1.
 \label{eq:circle-family}
\end{equation}
The experiments use $\lambda=-1/8,0,1/2,1$. The mean set is $\{\pm\arccos(11/12)\}$ at
$\lambda=-1/8$ and $\{0\}$ otherwise, with a quartic
minimum at $\lambda=0$ and positive curvature for
$\lambda>0$.
We compare with ordinary pointwise $\chi^2$-calibrated
empirical-likelihood score inversion
\citep{KurisuOtsu2026} and a single-centre Wald arc
(\suppcref{sec:supp-circle-model}). Neither comparator is calibrated for
simultaneous coverage of the entire Fr\'echet mean set.
Global coverage is $98.2$--$99.9\%$; Local coverage
is $95.8$--$99.1\%$ for nonnegative $\lambda$ and
$93.4$--$96.5\%$ with two means.
For positive $\lambda$, empirical likelihood gives
$93.5$--$94.5\%$ coverage. The Wald arc is near
nominal at $\lambda=1$, undercovers at $\lambda=1/2$,
and serves only as a formal full-set comparator
at nonpositive $\lambda$.

In all 12 circular settings, Local has higher observed whole-set
coverage than ordinary score-inverted empirical likelihood.

\putbib[ref]
\end{bibunit}
\switchtosupplementnumbering
\begin{bibunit}[apalike]
\makeatletter
\def\@extra@b@citeb{@supp}
\gdef\@extra@binfo{@supp}
\makeatother
\begin{center}
\begin{minipage}{0.94\textwidth}
\centering
{\fontsize{19}{23}\selectfont\bfseries Supplementary Material for\par
\vspace{2mm}
Inference for Fr\'echet Mean Sets\par}
\end{minipage}

\vspace{0.6cm}
{\fontsize{13.2}{16}\selectfont Ryunosuke Miyazaki\par}
\vspace{0.2cm}
{\fontsize{11}{13.2}\selectfont\itshape Department of Social Data Science, Hitotsubashi University\par}
\end{center}
\vspace{0.35cm}
\providecommand{\proofheading}[1]{Proof of #1}

\section{Preliminaries}
\subsection{Notation and loss geometry}
\label{sec:supp-loss}
Suppress the nominal level $\alpha$ from critical values and confidence sets.
For measurable $f$, write $Pf\coloneqq\int f\,dP$ and
$\|f\|_2\coloneqq(Pf^2)^{1/2}$. Define
$\mathcal I_P\coloneqq\{(q,r):\sigma_P(q,r)>0\}$ and
$\mathcal H_{P,n}^{\mathrm L}\coloneqq
\{h_{q,r}:(q,r)\in\wh{\mathcal B}_n^{\mathrm L}\cap\mathcal I_P\}$.
Fix $o\in\cQ$ and put
\[
 \Phi_P(q)\coloneqq d^2(\cdot,q)-d^2(\cdot,o)
                    -\{F_P(q)-F_P(o)\}\in L^2_0(P),
 \qquad \ell_{q,r}\coloneqq\Phi_P(q)-\Phi_P(r).
\]
Then $\sigma_P(q,r)=\|\ell_{q,r}\|_2$ and
$h_{q,r}=\ell_{q,r}/\sigma_P(q,r)$ on $\mathcal I_P$.
The second moment assumption gives, for $q,q'\in B(o,R)$,
\[
 \max\bigl(|F_P(q)-F_P(q')|,\|\Phi_P(q)-\Phi_P(q')\|_2\bigr)
 \leq2\|d(\cdot,o)+R\|_2\,d(q,q').
\]
Thus $F_P$, $\ell_{q,r}$ and $\sigma_P(q,r)$ are continuous, and
$(q,r)\mapsto h_{q,r}$ is $L^2(P)$-continuous on $\mathcal I_P$.

Write $\mathcal L(Z)$ for the law of $Z$; subscripts $P$ and $\xi$
indicate sampling and conditional multiplier laws; $\E_{P,\xi}$ averages
over sample and multipliers. For every $a,b>0$, the conventions are
\[
 \begin{aligned}
 R_n=o_{p_\xi}(1)\text{ in probability}
 &\iff\pr_P(\pr_\xi(|R_n|>a\mid Y_1,\ldots,Y_n)>b)\to0,\\
 R_{n,P}=o_{\mathcal P}(1)
 &\iff\sup_{P\in\mathcal P}\pr_P(|R_{n,P}|>a)\to0,\\
 R_{n,P}=o_{\mathcal P,\xi}(1)
 &\iff\sup_{P\in\mathcal P}
 \pr_P(\pr_\xi(|R_{n,P}|>a\mid Y_1,\ldots,Y_n)>b)\to0.
 \end{aligned}
\]
Here $O_{\mathcal P}(1)$ denotes uniform tightness, and
$d_{\mathrm{BL}}$ denotes bounded-Lipschitz distance between probability laws.
For laws $\mu,\nu$ on $\mathbb R$, write
$d_K(\mu,\nu)\coloneqq\sup_{t\in\mathbb R}
|\mu((-\infty,t])-\nu((-\infty,t])|$ (Kolmogorov distance).
For a metric $\delta$ and nonempty sets, write
$\vec d_\delta(A,B)\coloneqq\sup_{a\in A}\inf_{b\in B}\delta(a,b)$ and
$d_{H,\delta}(A,B)\coloneqq\max(\vec d_\delta(A,B),\vec d_\delta(B,A))$.
Set $\vec d_\delta(\varnothing,B)=d_{H,\delta}(\varnothing,\varnothing)=0$;
for nonempty $A$, set $\vec d_\delta(A,\varnothing)
=d_{H,\delta}(A,\varnothing)=d_{H,\delta}(\varnothing,A)=+\infty$.
$\operatorname{dist}$ and $d_H$ denote Euclidean point-to-set and Hausdorff distances.

\begin{lemma}[Simultaneous zero-variance comparisons]
\label{lem:zero-variance}
For every $(q,r)$ with $\sigma_P(q,r)=0$,
$d^2(y,q)-d^2(y,r)=F_P(q)-F_P(r)$ at every
$y\in\operatorname{supp}(P)$. Consequently, a.s.\ and
simultaneously for all such pairs,
$F_n(q)-F_n(r)=F_P(q)-F_P(r)$ and
$\wh\sigma_n(q,r)=0$.
\end{lemma}

\begin{proof}
On a Polish space, $P(\operatorname{supp}(P))=1$
\citep[Lemma~1.19]{Kallenberg2021}. For any zero-variance pair,
$\ell_{q,r}$ is continuous in the observation and vanishes $P$-a.s.; its open nonzero set is therefore disjoint from the support.
On the single probability-one event that all observations lie in this
support, the empirical identities hold for every such pair.
\end{proof}

\begin{lemma}[Squared-distance contrast scales]
\label{lem:automatic-linear-scale}
For any probability measure $R$ with finite second moment, define
$F_R(q)=\int d^2(y,q)\,dR(y)$ and
$\sigma_R(q,r)^2=\Var_R(d^2(Y,q)-d^2(Y,r))$. Then
\begin{equation}
 \begin{aligned}
 \sigma_R(q,r)
 &\le2d(q,r)\min(F_R(q)^{1/2},F_R(r)^{1/2})\\
 &\le d(q,r)\{F_R(q)^{1/2}+F_R(r)^{1/2}\}
 \le2\{F_R(q)+F_R(r)\}.
 \end{aligned}
 \label{eq:loss-risk-scale}
\end{equation}
This applies to $R=P$ and to the empirical law $P_n$.
For $m\in M_P$ and $(m,r)\in\mathcal I_P$,
\begin{equation}
 F_P(r)-F_P^\star\le4\lambda_P(m,r)\{1+\lambda_P(m,r)\}F_P^\star.
 \label{eq:small-slack-risk}
\end{equation}
For every $q\in\cQ$,
\[
 \inf_{m\in M_P}\sigma_P(q,m)
 \le2(F_P^\star)^{1/2}d(q,M_P).
\]
\end{lemma}
\begin{proof}
The triangle and reverse triangle inequalities give
\[
 |d(y,q)^2-d(y,r)^2-d(q,r)^2|\le2d(q,r)d(y,r).
\]
Variance is bounded by the mean square about any constant. Interchanging
$q,r$ gives the first bound in \eqref{eq:loss-risk-scale}; use
$d(q,r)\le F_R(q)^{1/2}+F_R(r)^{1/2}$ and
$(a+b)^2\le2(a^2+b^2)$ for the remaining bounds.
For \eqref{eq:small-slack-risk},
\[
 \begin{aligned}
 \{F_P(r)^{1/2}-(F_P^\star)^{1/2}\}
 \{F_P(r)^{1/2}+(F_P^\star)^{1/2}\}
 &=\lambda_P(m,r)\sigma_P(m,r)\\
 &\le2\lambda_P(m,r)(F_P^\star)^{1/2}
       \{F_P(r)^{1/2}+(F_P^\star)^{1/2}\}.
 \end{aligned}
\]
The final factor is positive on $\mathcal I_P$; cancellation and
squaring yield the claim. Taking the infimum of the first
bound in \eqref{eq:loss-risk-scale} over $m\in M_P$ proves the final assertion.
\end{proof}

For $\eta>0$, put
\[
 \mathcal A_P(\eta)\coloneqq
 \{h_{m,r}:m\in M_P,\ (m,r)\in\mathcal I_P,\ \lambda_P(m,r)\leq\eta\},
 \qquad
 \mathcal H_{P,\eta}\coloneqq\overline{\mathcal A_P(\eta)}^{\,L^2(P)}.
\]
Thus $\mathcal H_{P,0}=\bigcap_{\eta>0}\mathcal H_{P,\eta}$.
Put $M_n\coloneqq\argmin_{q\in\cQ}F_n(q)$, choosing
$\widehat m_n\in M_n$ only when $M_n\ne\varnothing$.
If $F_n^\star>0$, the infimum definition gives $\wh M_n\ne\varnothing$;
if $F_n^\star=0$, all observations coincide and $\wh M_n=\{Y_1\}$.
For generic $f$, put
\[
 \mathbb G_nf\coloneqq n^{-1/2}\sum_{i=1}^n\{f(Y_i)-Pf\},\qquad
 \mathbb G_n^\#f\coloneqq n^{-1/2}\sum_{i=1}^n
 \xi_i\left\{f(Y_i)-\frac1n\sum_{j=1}^n f(Y_j)\right\}.
\]
If $\mathcal I_P=\varnothing$, \cref{lem:zero-variance} gives
$\mathcal C_n(c)=M_P$ for every $c\ge0$ a.s..

Scalar suprema and selected-set events use completed probability spaces:
continuity of the loss maps on the Polish parameter space and Borel
selector graphs give universally measurable versions. Process-level
probabilities and expectations are outer when needed; bounded-Lipschitz
distances use the expectation formulation, conditional over multipliers,
as in \citet[Sections~1.12 and~2.9]{VaartWellner2023}.
Empirical and multiplier suprema use actual directions in $\mathcal H_P$;
only the Gaussian process is extended to
$\overline{\mathcal H_P}^{\,L^2(P)}$.

\subsection{Uniform studentization}
\label{sec:supp-studentization}
Fix a nonempty distribution class $\mathcal P$ satisfying
\maincref{ass:main-process}.
All uniform bounds, limits and stochastic orders in this subsection
refer to this class.
Let $\mathrm{BL}_1(P)$ be the maps from $\ell^\infty(\mathcal H_P)$
to $[-1,1]$ with Lipschitz constant at most one. The uniform Donsker
condition in \maincref{ass:main-process} means
\begin{equation}
 \sup_{P\in\mathcal P}\sup_{f\in\mathrm{BL}_1(P)}
 \left|\E_P[f(\mathbb G_n)]-\E[f(\mathbb G_P)]\right|\longrightarrow0.
 \label{eq:uniform-donsker}
\end{equation}
Uniform pre-Gaussianity means that each $\mathbb G_P$ has bounded,
uniformly $\rho_P$-continuous paths and
\begin{equation}
 \sup_{P\in\mathcal P}\E[\|\mathbb G_P\|_{\mathcal H_P}]<\infty,
 \qquad
 \lim_{\delta\downarrow0}\sup_{P\in\mathcal P}
 \E\left[\sup_{\substack{h,g\in\mathcal H_P\\\rho_P(h,g)<\delta}}
       |\mathbb G_P(h)-\mathbb G_P(g)|\right]=0.
 \label{eq:uniform-pregaussian}
\end{equation}
These conditions use the framework of
\citet[Definitions~2.0--2.1]{SheehyWellner1992} and
\citet[Section~2.8.2]{VaartWellner2023}, with finite nets in each
$L^2(P)$ and bounds uniform over $P$.

\begin{lemma}[Uniform process approximation and second moments]
\label{lem:uniform-squares}
Under \maincref{ass:main-process},
the following hold uniformly over $P\in\mathcal P$.
\begin{enumerate}
\item[(i)] The empirical processes satisfy
\[
 \|\mathbb G_n\|_{\mathcal H_P}=O_{\mathcal P}(1),
\]
and, for every $a>0$,
\[
 \lim_{\delta\downarrow0}\limsup_{n\to\infty}\sup_{P\in\mathcal P}
 \pr_P\left(\sup_{\substack{h,g\in\mathcal H_P\\\rho_P(h,g)<\delta}}
                 |\mathbb G_n(h-g)|>a\right)=0.
\]
\item[(ii)] The conditional multiplier laws satisfy
\begin{equation}
 d_{\mathrm{BL}}(\mathcal L_\xi(\mathbb G_n^\#\mid Y_1,\ldots,Y_n),
                       \mathcal L(\mathbb G_P))
       =o_{\mathcal P}(1).
 \label{eq:uniform-multiplier-squares}
\end{equation}
In particular, the multiplier processes are conditionally uniformly tight and
asymptotically $\rho_P$-equicontinuous in probability, uniformly in $P$.
\item[(iii)] The empirical second moments satisfy
\[
 \sup_{h\in\mathcal H_P}|P_nh^2-1|=o_{\mathcal P}(1).
\]
\end{enumerate}
\end{lemma}
\begin{proof}
By \eqrefs{eq:uniform-donsker,eq:uniform-pregaussian}, the bounded-Lipschitz
modulus argument of \citet[Theorem~2.8.2]{VaartWellner2023} gives (i).
Sudakov's minoration and the first bound in
\eqref{eq:uniform-pregaussian} also give
\begin{equation}
 \sup_{P\in\mathcal P}N(\delta,\mathcal H_P,\rho_P)<\infty
 \qquad(\delta>0).
 \label{eq:uniform-empirical-control}
\end{equation}
These arguments apply to the $P$-dependent classes with bounds uniform in $P$.
For every $a>0$, the envelope satisfies
\[
 \sup_P\E_P\left[n^{-1}\max_{i\le n}H_P(Y_i)^2\right]
 \le a^2+\sup_P \E_P[H_P(Y)^2\mathds1\{H_P(Y)>an^{1/2}\}].
\]
Letting $n\to\infty$, then $a\downarrow0$, Hoffmann--J\o rgensen's
inequality \citep[Lemma~2.3.11 and Proposition~A.1.5]{VaartWellner2023}
and (i) give
\begin{equation}
 \begin{gathered}
 \sup_{n,P}\E_P[\|\mathbb G_n\|_{\mathcal H_P}]<\infty,\\
 \lim_{\delta\downarrow0}\limsup_n\sup_P
 \E_P\left[\sup_{\substack{h,g\in\mathcal H_P\\\rho_P(h,g)<\delta}}
                |\mathbb G_n(h-g)|\right]=0.
 \end{gathered}
 \label{eq:uniform-empirical-means}
\end{equation}

For $1\le n_0\le n$, the multiplier inequality
\citep[Lemma~2.9.1]{VaartWellner2023}, followed by symmetrization, yields
for a universal constant $C$
\begin{equation}
\begin{split}
 &\sup_P\E_{P,\xi}\left[\sup_{\substack{h,g\in\mathcal H_P\\\rho_P(h,g)<\delta}}
 \left|\frac1{n^{1/2}}\sum_{i=1}^n\xi_i(h-g)(Y_i)\right|\right]\\
 &\quad\le
 \frac{C(n_0-1)}{n^{1/2}}\sup_P PH_P\,
              \E[\max_{i\le n}|\xi_i|]
 +C\sup_{k\ge n_0}\sup_P
       \E_P\left[\sup_{\substack{h,g\in\mathcal H_P\\\rho_P(h,g)<\delta}}
                  |\mathbb G_k(h-g)|\right].
\end{split}
 \label{eq:uniform-multiplier-modulus}
\end{equation}
Since $\E[\max_{i\le n}|\xi_i|]=O(\{\log(n+1)\}^{1/2})$,
let $n\to\infty$, $n_0\to\infty$, then $\delta\downarrow0$.
Equation \eqref{eq:uniform-empirical-means} and Markov's inequality
give conditional equicontinuity of the uncentred multiplier process
$\widetilde{\mathbb G}_n^\#h\coloneqq
n^{-1/2}\sum_i\xi_i h(Y_i)$, uniformly in probability.
Since $Ph=0$ for every $h\in\mathcal H_P$, empirical centring adds only
\[
 \|\mathbb G_n^\#-\widetilde{\mathbb G}_n^\#\|_{\mathcal H_P}
 =\frac{\|\mathbb G_n\|_{\mathcal H_P}}{n^{1/2}}
   \left|\frac1{n^{1/2}}\sum_i\xi_i\right|
 =o_{\mathcal P,\xi}(1).
\]
Thus the same conditional equicontinuity holds for $\mathbb G_n^\#$.

For fixed $\delta>0$, choose a $\rho_P$-net
$h_1,\ldots,h_N\in\mathcal H_P$ with
$N\le N_\delta<\infty$ uniformly in $P$, using
\eqref{eq:uniform-empirical-control}.
Truncation, Chebyshev's inequality and a union bound give, for $L,\epsilon>0$,
\begin{equation}
 \begin{split}
 &\sup_P\pr_P\left(
 \max_{j,k\le N}\left|
 \frac1n\sum_i h_j(Y_i)h_k(Y_i)-P(h_jh_k)
 \right|>\epsilon\right)\\
 &\quad\le
 \frac{CN_\delta^2L^4}{n\epsilon^2}
 +\frac{CN_\delta^2}{\epsilon}\sup_P \E_P[H_P(Y)^2\mathds1\{H_P(Y)>L\}].
 \end{split}
 \label{eq:uniform-net-covariance}
\end{equation}
Letting $n\to\infty$, then $L\to\infty$, gives uniform covariance
consistency on these nets. The centring correction
$(P_nh_j)(P_nh_k)$ is uniformly $o_{\mathcal P}(1)$ by (i).
Conditionally on the observations, the multiplier vector on each net
is centred Gaussian. Its target covariance matrices have diagonal entries
one and form a bounded subset of the positive semidefinite matrices in
dimension at most $N_\delta$. Continuity of the matrix square root is
uniform on compact subsets, including at singular matrices; coupling
the two Gaussian vectors with a common standard normal vector therefore
gives uniform conditional bounded-Lipschitz convergence on the nets.

Project each $h\in\mathcal H_P$ onto a net point within $\delta$.
For any map in $\mathrm{BL}_1(P)$, the resulting approximation error
is bounded by the conditional expectation of the multiplier modulus,
truncated at two; the Gaussian error has the analogous bound.
Conditional equicontinuity makes the former vanish uniformly in probability
as $n\to\infty$, then $\delta\downarrow0$, and
\eqref{eq:uniform-pregaussian} controls the latter.
This uniform version of the finite-net argument in
\citet[Theorem~2.9.6]{VaartWellner2023} proves
\eqref{eq:uniform-multiplier-squares}. Conditional uniform tightness
follows from this approximation and \eqref{eq:uniform-pregaussian}, proving (ii).

Conditionally on $Y_1,\ldots,Y_n$, $\mathbb G_n^\#h$ is centred Gaussian
with variance $V_{n,P}(h)\coloneqq P_nh^2-(P_nh)^2$.
Since $Ph^2=1$, the conditional characteristic function and the
$\mathrm{BL}_1(P)$ map $z\mapsto\cos\{z(h)\}$ give, by (ii),
\[
 \E_\xi[\cos(\mathbb G_n^\#h)\mid Y_1,\ldots,Y_n]=e^{-V_{n,P}(h)/2},
 \qquad
 \sup_h|e^{-V_{n,P}(h)/2}-e^{-1/2}|=o_{\mathcal P}(1).
\]
Continuity of the inverse map $x\mapsto-2\log x$ near $e^{-1/2}$ gives
$\sup_h|V_{n,P}(h)-1|=o_{\mathcal P}(1)$, whence
\[
 \sup_{h\in\mathcal H_P}|P_nh^2-1|
 \leq\sup_{h\in\mathcal H_P}|V_{n,P}(h)-1|+\frac{\|\mathbb G_n\|_{\mathcal H_P}^2}{n}
 =o_{\mathcal P}(1).
\]
This proves (iii).
\end{proof}

\begin{lemma}[Uniform studentization]
\label{lem:relative-variance}
Under \maincref{ass:main-process},
\begin{equation*}
 \sup_{(q,r)\in\mathcal I_P}
 \left|\frac{\wh\sigma_n(q,r)}{\sigma_P(q,r)}-1\right|
       =o_{\mathcal P}(1),
 \qquad
 \inf_{P\in\mathcal P}
 \pr_P(\{(q,r):\wh\sigma_n(q,r)>0\}=\mathcal I_P)\longrightarrow1.
\end{equation*}
\end{lemma}
\begin{proof}
For every $(q,r)\in\mathcal I_P$,
\begin{equation}
 \frac{\wh\sigma_n^2(q,r)}{\sigma_P^2(q,r)}
 =\frac1n\sum_i h_{q,r}(Y_i)^2
       -\left(\frac1n\sum_i h_{q,r}(Y_i)\right)^2
 =V_{n,P}(h_{q,r}).
 \label{eq:proof-relative-variance}
\end{equation}
Apply \cref{lem:uniform-squares}(i),(iii) and take square roots.
On $\{\sup_{\mathcal I_P}|\wh\sigma_n/\sigma_P-1|<1/2\}$,
\cref{lem:zero-variance} gives equality of the two index sets.
\end{proof}
\subsection{Suprema and quantiles}
\label{sec:supp-suprema}
\begin{lemma}[Suprema over estimated index sets]
\label{lem:moving-sets}
Fix $P$ and assume \maincref{ass:main-process} for $\mathcal P=\{P\}$.
All random sets below are measurable with respect to $Y_1,\ldots,Y_n$.
\begin{enumerate}
\item[(i)] Let $A\subseteq\overline{\mathcal H_P}^{\,L^2(P)}$ be fixed, nonempty and
compact in $\rho_P$. If $A_n\subseteq\mathcal H_P$ satisfies
\[
 d_{H,\rho_P}(\overline{A_n}^{\,L^2(P)},A)\to_p0,
\]
then
\[
 \sup_{h\in A_n}\mathbb G_nh\rightsquigarrow
 \sup_{h\in A}\mathbb G_P(h)
\]
and
\[
 d_{\mathrm{BL}}\left(
 \mathcal L_\xi\left(\sup_{h\in A_n}\mathbb G_n^\#h\mid Y_1,\ldots,Y_n\right),
 \mathcal L\left(\sup_{h\in A}\mathbb G_P(h)\right)\right)\to_p0.
\]
\item[(ii)] If $A_n,B_n\subseteq\mathcal H_P$ are nonempty with
probability tending to one and
\[
 d_{H,\rho_P}(\overline{A_n}^{\,L^2(P)},
             \overline{B_n}^{\,L^2(P)})\to_p0,
\]
then
\[
 \sup_{h\in A_n}\mathbb G_nh-\sup_{h\in B_n}\mathbb G_nh=o_p(1),
\]
and
\[
 \sup_{h\in A_n}\mathbb G_n^\#h-\sup_{h\in B_n}\mathbb G_n^\#h
 =o_{p_\xi}(1)\quad\text{in probability}.
\]
\end{enumerate}
\end{lemma}
\begin{proof}
Assign zero to suprema over empty random sets; under the assumptions
of each part, these events have probability tending to zero.
For (i), fix $\delta>0$ and choose a finite $\delta$-net
$a_1,\ldots,a_N\in A$. Since $A\subseteq\overline{\mathcal H_P}^{\,L^2(P)}$, choose
$g_j\in\mathcal H_P$ with $\rho_P(g_j,a_j)<\delta$.
On $\{d_{H,\rho_P}(\overline{A_n}^{\,L^2(P)},A)<\delta\}$,
every $h\in A_n$ is within $3\delta$ of some $g_j$, while every
$g_j$ is within $2\delta$ of some $h\in A_n$.
Here distance to $A_n$ equals distance to its $L^2(P)$ closure,
so both comparisons use actual elements of $\mathcal H_P$. Hence
\begin{equation}
 \left|\sup_{h\in A_n}\mathbb G_nh-\max_{j\le N}\mathbb G_ng_j\right|
 \le
 \sup_{\substack{h,g\in\mathcal H_P\\\rho_P(h,g)<4\delta}}
                      |\mathbb G_n(h-g)|.
 \label{eq:empirical-modulus}
\end{equation}
The same inequality holds with $\mathbb G_n$ replaced by $\mathbb G_n^\#$.
The Gaussian extension to $\overline{\mathcal H_P}^{\,L^2(P)}$ satisfies
\[
 \left|\sup_{h\in A}\mathbb G_P(h)-\max_{j\le N}\mathbb G_P(g_j)\right|
 \le
 \sup_{\substack{h,g\in\overline{\mathcal H_P}^{\,L^2(P)}\\\rho_P(h,g)<2\delta}}
                      |\mathbb G_P(h)-\mathbb G_P(g)|.
\]
For fixed $\delta$, the empirical and conditional multiplier laws of
the finite maxima converge to their Gaussian counterparts by
\eqref{eq:uniform-donsker} and \cref{lem:uniform-squares}(ii), since
$z\mapsto\max_{j\le N}z(g_j)$ is $1$-Lipschitz in the supremum norm.
The empirical and conditional multiplier equicontinuity in
\cref{lem:uniform-squares}(i),(ii), together with Gaussian path continuity,
make the approximation errors vanish as $n\to\infty$, then
$\delta\downarrow0$. The exceptional event is determined by the data
and has probability tending to zero. This proves (i); the same
Hausdorff condition implies $\pr_P(A_n=\varnothing)\to0$.

For (ii), on the event that both sets are nonempty and their closures
have Hausdorff distance less than $\delta$, every point of either set
is within $2\delta$ of a point of the other. Consequently,
\[
 \left|\sup_{h\in A_n}\mathbb G_nh-\sup_{g\in B_n}\mathbb G_ng\right|
 \le
 \sup_{\substack{h,g\in\mathcal H_P\\\rho_P(h,g)<2\delta}}
                     |\mathbb G_n(h-g)|.
\]
The same bound holds for $\mathbb G_n^\#$.
The exceptional event is determined by the data and has probability
tending to zero. Empirical and conditional multiplier equicontinuity
from \cref{lem:uniform-squares}(i),(ii) prove both assertions by letting
$n\to\infty$, then $\delta\downarrow0$.
\end{proof}

\begin{lemma}[Gaussian supremum quantiles]
\label{lem:gaussian-quantile}
Fix $P$ and assume \maincref{ass:main-process} for $\mathcal P=\{P\}$.
For each fixed nonempty $\rho_P$-compact set
$A\subseteq\overline{\mathcal H_P}^{\,L^2(P)}$, put
$Z_A\coloneqq\sup_{h\in A}\mathbb G_P(h)$.
For $1/2<p<1$, write $c_{A,p}\coloneqq q_p(Z_A)$ and
$z_p\coloneqq\Phi^{-1}(p)$, where $\Phi$ is the standard normal
distribution function.
\begin{enumerate}
\item[(i)] The law of $Z_A$ is atomless, its distribution function is
strictly increasing on $(0,\infty)$, and $\pr(Z_A\le0)\le1/2$.
Moreover, $c_{A,p}\ge z_p>0$ and $\pr(Z_A\le c_{A,p})=p$.
For every $\epsilon>0$, the following bounds hold independently of $P$ and $A$:
\begin{equation}
 \begin{aligned}
 \pr(Z_A\le c_{A,p}-\epsilon)&\le\Phi(z_p-\epsilon)<p,\\
 \pr(Z_A\le c_{A,p}+\epsilon)&\ge\Phi(z_p+\epsilon)>p.
 \end{aligned}
 \label{eq:uniform-quantile-separation}
\end{equation}
\item[(ii)] If random probability laws $\nu_n$ on $\mathbb R$ converge in
bounded-Lipschitz distance in probability to either $\mathcal L(Z_A)$ or
$\mathcal L((Z_A)_+)$, then $q_p(\nu_n)\to_p c_{A,p}$.
\item[(iii)] If $A\subsetneq B\subseteq
\overline{\mathcal H_P}^{\,L^2(P)}$ and $B$ is $\rho_P$-compact, then
$\pr(Z_B\le t)<\pr(Z_A\le t)$ for every $t>0$ and
$c_{A,p}<c_{B,p}$ for $1/2<p<1$.
\end{enumerate}
\end{lemma}
\begin{proof}
Every index has $P$-mean zero and $L^2(P)$ norm one, and
\eqref{eq:uniform-pregaussian} gives bounded continuous Gaussian paths on $A$.
Gaussian anti-concentration
\citep[Theorem~2.1]{ChernozhukovChetverikovKato2014} therefore makes $Z_A$
atomless. As a centred Gaussian random element of $C(A)$, the process has
positive probability in every open ball about zero
\citep[Lemma~5.1]{VaartZanten2008}. Thus, for $h_0\in A$ and $t>0$,
\[
 0<\pr\left(\sup_{h\in A}|\mathbb G_P(h)|<t\right)
 \le\pr(Z_A\le t)\le\pr(\mathbb G_P(h_0)\le t)=\Phi(t)<1.
\]
The same coordinate comparison gives $\pr(Z_A\le0)\le1/2$.

For finite $A$, its Gaussian supremum is a maximum of unit linear
projections of a standard Gaussian vector, hence is $1$-Lipschitz.
Gaussian isoperimetry applied to its sublevel sets \citep{Borell1975} gives
\[
 \pr(Z_A\le t+u)
 \ge\Phi\left(\Phi^{-1}\{\pr(Z_A\le t)\}+u\right),
 \qquad t\in\mathbb R,\quad u>0,
\]
where $\Phi^{-1}(0)=-\infty$. Increasing finite subsets with dense union
and path continuity extend this bound to compact $A$.
It implies strict increase on $(0,\infty)$; atomlessness and the coordinate
comparison give $\pr(Z_A\le c_{A,p})=p$ and $c_{A,p}\ge z_p>0$.
Applying the bound at $t=c_{A,p}$ and $t=c_{A,p}-\epsilon$, with
increment $\epsilon$, proves \eqref{eq:uniform-quantile-separation} and (i).

For (ii), fix $0<\epsilon<c_{A,p}$. The laws of $Z_A$ and $(Z_A)_+$ have
the same distribution function at positive arguments, and both are
continuous at $c_{A,p}\pm\epsilon$. Bounded-Lipschitz convergence and
\eqref{eq:uniform-quantile-separation} therefore give, with probability
tending to one,
\[
 \nu_n((-\infty,c_{A,p}-\epsilon])<p
 <\nu_n((-\infty,c_{A,p}+\epsilon]),
\]
and hence $|q_p(\nu_n)-c_{A,p}|\le\epsilon$.

For (iii), fix $t>0$ and choose $h_*\in B\setminus A$.
Compactness and unit norms give
$\sup_{h\in A}\langle h,h_*\rangle_{L^2(P)}<1$.
Choose $s>t$ with $s\sup_{h\in A}\langle h,h_*\rangle_{L^2(P)}<t$ and define
\[
 R(h)\coloneqq\mathbb G_P(h)
 -\langle h,h_*\rangle_{L^2(P)}\mathbb G_P(h_*),\qquad h\in A.
\]
This centred continuous Gaussian process is independent of
$\mathbb G_P(h_*)\sim N(0,1)$, since all cross-covariances vanish.
Choose $\epsilon>0$ with $s-\epsilon>t$ and
$s\sup_{h\in A}\langle h,h_*\rangle_{L^2(P)}+2\epsilon<t$.
Independence and small-ball positivity give
\[
 \pr\left(|\mathbb G_P(h_*)-s|<\epsilon,\quad
          \sup_{h\in A}|R(h)|<\epsilon\right)>0.
\]
On this event, using $|\langle h,h_*\rangle_{L^2(P)}|\le1$,
\[
 Z_A\le s\sup_{h\in A}\langle h,h_*\rangle_{L^2(P)}+2\epsilon
 <t<\mathbb G_P(h_*)\le Z_B.
\]
Since $Z_A\le Z_B$ a.s.,
\[
 \pr(Z_A\le t)-\pr(Z_B\le t)=\pr(Z_A\le t<Z_B)>0.
\]
At $t=c_{A,p}>0$, continuity of the distribution function of $Z_B$
therefore yields $c_{B,p}>c_{A,p}$.
\end{proof}

\begin{lemma}[Uniform comparison of supremum quantiles]
\label{lem:uniform-quantiles}
Fix a nonempty class $\mathcal P$ satisfying
\maincref{ass:main-process}. For arbitrary nonempty deterministic
pair sets $\mathcal D_{n,P}\subseteq\mathcal I_P$, possibly depending on
$n$ and $P$, put $A_{n,P}=\{h_{q,r}:(q,r)\in\mathcal D_{n,P}\}$ and
$Z_{n,P}=\sup_{h\in A_{n,P}}\mathbb G_P(h)$.
On the event that the empirical and population variance index sets agree and
\[
 \frac12\le\frac{\wh\sigma_n(q,r)}{\sigma_P(q,r)}\le2
 \qquad((q,r)\in\mathcal I_P),
\]
set
\begin{equation*}
 U_{n,P}=\sup_{(q,r)\in\mathcal D_{n,P}}
       \frac{\sigma_P(q,r)}{\wh\sigma_n(q,r)}\mathbb G_nh_{q,r},
 \quad
 V_{n,P}^\#=\sup_{(q,r)\in\mathcal D_{n,P}}
       \frac{\sigma_P(q,r)}{\wh\sigma_n(q,r)}\mathbb G_n^\#h_{q,r},
\end{equation*}
and assign both quantities zero on its complement.
The following hold uniformly over these choices.
\begin{enumerate}
\item[(i)] The empirical and conditional multiplier laws satisfy
\begin{equation*}
 \sup_P d_K(\mathcal L_P(U_{n,P}),\mathcal L(Z_{n,P}))\to0,
 \qquad
 d_K(\mathcal L_\xi(V_{n,P}^\#\mid Y_1,\ldots,Y_n),
                  \mathcal L(Z_{n,P}))=o_{\mathcal P}(1).
\end{equation*}
\item[(ii)] For $0<\alpha<1/2$, let $q_{n,P}^\#$ be the conditional
$(1-\alpha)$-quantile of $(V_{n,P}^\#)_+$. Then
\begin{equation*}
 \liminf_n\inf_P\pr_P((U_{n,P})_+\le q_{n,P}^\#)
       \ge1-\alpha.
\end{equation*}
\end{enumerate}
\end{lemma}
\begin{proof}
For $a>0$, Gaussian anti-concentration
\citep[Theorem~2.1]{ChernozhukovChetverikovKato2014} gives
\begin{equation}
 \sup_{P,n,t}\pr(|Z_{n,P}-t|\le a)
 \le4a\left(\sup_P\E[\|\mathbb G_P\|_{\mathcal H_P}]+1\right).
 \label{eq:uniform-anticoncentration}
\end{equation}
On the defining event,
\cref{lem:uniform-squares,lem:relative-variance} give
\[
 \begin{split}
 \left|U_{n,P}-\sup_{h\in A_{n,P}}\mathbb G_nh\right|
 &\le \|\mathbb G_n\|_{\mathcal H_P}
 \sup_{\mathcal I_P}\left|\frac{\sigma_P}{\wh\sigma_n}-1\right|
 =o_{\mathcal P}(1),\\
 \left|V_{n,P}^\#-\sup_{h\in A_{n,P}}\mathbb G_n^\#h\right|
 &\le \|\mathbb G_n^\#\|_{\mathcal H_P}
 \sup_{\mathcal I_P}\left|\frac{\sigma_P}{\wh\sigma_n}-1\right|
 =o_{\mathcal P,\xi}(1).
 \end{split}
\]
The complement is uniformly negligible. The supremum map over each
deterministic $A_{n,P}$ is 1-Lipschitz, so \cref{lem:uniform-squares} and
these bounds give bounded-Lipschitz approximation of both laws in (i).
For any law $\nu$ on $\mathbb R$ and $0<a\le1$, smoothing and
\eqref{eq:uniform-anticoncentration} give, for a universal constant $C$,
\[
 d_K(\nu,\mathcal L(Z_{n,P}))
 \le \frac{C}{a}d_{\mathrm{BL}}(\nu,\mathcal L(Z_{n,P}))
   +4a\left(\sup_P\E[\|\mathbb G_P\|_{\mathcal H_P}]+1\right).
\]
Letting $n\to\infty$ and then $a\downarrow0$ proves (i).

For (ii), fix $0<\gamma<1/2-\alpha$.
The closure $\overline{A_{n,P}}^{\,L^2(P)}$ is compact by
\eqref{eq:uniform-empirical-control}, so \cref{lem:gaussian-quantile} gives
\[
 c_{n,P}^-\coloneqq q_{1-\alpha-\gamma}(Z_{n,P})
 \ge\Phi^{-1}(1-\alpha-\gamma)>0.
\]
Let $E_{n,P}$ be the event that the conditional distance in (i) is at
most $\gamma$. On $E_{n,P}$, the clipped conditional distribution function
is below $1-\alpha$ for every $t<c_{n,P}^-$, so $q_{n,P}^\#\ge c_{n,P}^-$.
Since $c_{n,P}^->0$, (i) and \cref{lem:gaussian-quantile} give
\[
 \pr_P((U_{n,P})_+\le q_{n,P}^\#)
 \ge\pr_P(U_{n,P}\le c_{n,P}^-)-\pr_P(E_{n,P}^c)
 \ge1-\alpha-\gamma-o(1),
\]
uniformly over $P$ and the deterministic pair sets.
Letting $\gamma\downarrow0$ proves (ii).
\end{proof}

\begin{lemma}[Compact binding localization]
\label{lem:compact-binding}
Suppose $\overline{\mathcal H_P}^{\,L^2(P)}$ is compact.
If $\mathcal H_{P,0}\neq\varnothing$,
then
\begin{equation*}
 d_{H,\rho_P}(\mathcal H_{P,\eta},\mathcal H_{P,0})\longrightarrow0
 \quad(\eta\downarrow0).
\end{equation*}
If $\mathcal H_{P,0}=\varnothing$, then $\mathcal A_P(\eta_0)$ is empty
for some $\eta_0>0$.
\end{lemma}
\begin{proof}
The sets $\mathcal H_{P,\eta}$ are decreasing closed subsets of the
compact space $\overline{\mathcal H_P}^{\,L^2(P)}$, hence compact.

If $\mathcal H_{P,0}\ne\varnothing$ and the Hausdorff convergence failed,
there would exist $\eta_j\downarrow0$, $\varepsilon>0$, and
$h_j\in\mathcal H_{P,\eta_j}$ such that
\[
 \inf_{h\in\mathcal H_{P,0}}\rho_P(h_j,h)\ge\varepsilon.
\]
Compactness gives a convergent subsequence with limit
$h\in\overline{\mathcal H_P}^{\,L^2(P)}$.
For every fixed $\eta>0$, eventually $h_j\in\mathcal H_{P,\eta}$;
closedness therefore gives $h\in\mathcal H_{P,\eta}$.
Thus $h\in\mathcal H_{P,0}$, a contradiction.

If $\mathcal H_{P,0}=\varnothing$ and every $\mathcal H_{P,1/j}$ were
nonempty, the finite-intersection property of compact sets would imply
\[
 \varnothing\ne\bigcap_{j\ge1}\mathcal H_{P,1/j}=\mathcal H_{P,0},
\]
again a contradiction. Hence $\mathcal H_{P,\eta_0}=\varnothing$ for some
$\eta_0>0$, and therefore $\mathcal A_P(\eta_0)=\varnothing$.
\end{proof}

\section{Proofs of the main results}
\subsection{Binding Gaussian limits}
\begin{proof}[\proofheading{\mainCref{thm:binding-main}}]
Apply the process and studentization lemmas with $\mathcal P=\{P\}$.
By \cref{lem:relative-variance}, with probability tending to one,
the empirical and population variance index sets agree and
\[
 \frac12\le\frac{\sigma_P(q,r)}{\wh\sigma_n(q,r)}\le2
 \qquad ((q,r)\in\mathcal I_P).
\]
On this event, put
\[
 \varepsilon_n\coloneqq\sup_{(q,r)\in\mathcal I_P}
 \left|\frac{\sigma_P(q,r)}{\wh\sigma_n(q,r)}-1\right|,
\]
and set $\varepsilon_n=0$ on its complement. Then
$\varepsilon_n=o_p(1)$. Equation
\mainequationref{eq:studentized-mean-decomposition} and
$|\mathbb G_nh_{m,r}|\le \|\mathbb G_n\|_{\mathcal H_P}$ give, for $m\in M_P$ and
$(m,r)\in\mathcal I_P$,
\[
 \mathbb G_nh_{m,r}-\|\mathbb G_n\|_{\mathcal H_P}\varepsilon_n
       -2n^{1/2}\lambda_P(m,r)
 \le T_n(m,r)
 \le\mathbb G_nh_{m,r}+\|\mathbb G_n\|_{\mathcal H_P}\varepsilon_n
       -\tfrac12n^{1/2}\lambda_P(m,r).
\]
If $\mathcal H_{P,0}\ne\varnothing$, restricting the lower bound to
$\lambda_P(m,r)\le n^{-1}$ and taking the supremum yields
\begin{equation}
 \sup_{h\in\mathcal A_P(n^{-1})}\mathbb G_nh
 -\|\mathbb G_n\|_{\mathcal H_P}\varepsilon_n-2n^{-1/2}
 \le T_n(M_P).
 \label{eq:binding-lower}
\end{equation}
For the upper bound, drop the drift when
$\lambda_P(m,r)\le n^{-1/4}$; otherwise,
\[
 T_n(m,r)\le2\|\mathbb G_n\|_{\mathcal H_P}-\tfrac12n^{1/2}\lambda_P(m,r)
          \le2\|\mathbb G_n\|_{\mathcal H_P}-\tfrac12n^{1/4}.
\]
Hence
\begin{equation}
 T_n(M_P)\le\max\left(
 \sup_{h\in\mathcal A_P(n^{-1/4})}\mathbb G_nh
 +\|\mathbb G_n\|_{\mathcal H_P}\varepsilon_n,
 \ 2\|\mathbb G_n\|_{\mathcal H_P}-\tfrac12n^{1/4}\right).
 \label{eq:binding-upper}
\end{equation}
By \cref{lem:compact-binding,lem:moving-sets},
\[
 \sup_{h\in\mathcal A_P(n^{-1})}\mathbb G_nh\rightsquigarrow T_P^0,
 \qquad
 \sup_{h\in\mathcal A_P(n^{-1/4})}\mathbb G_nh
 -\sup_{h\in\mathcal A_P(n^{-1})}\mathbb G_nh=o_p(1).
\]
Since $\|\mathbb G_n\|_{\mathcal H_P}\varepsilon_n=o_p(1)$,
\eqrefs{eq:binding-lower,eq:binding-upper} prove the limit.

If $\mathcal H_{P,0}=\varnothing$, \cref{lem:compact-binding} gives
$\eta_0>0$ with $\mathcal A_P(\eta_0)=\varnothing$. On the same event,
\[
 T_n(M_P)\le2\|\mathbb G_n\|_{\mathcal H_P}-\tfrac12n^{1/2}\,\eta_0\to_p-\infty.
\]
\end{proof}

\begin{lemma}[Uniform approximation on binding directions]
\label{lem:uniform-binding-transfer}
Let $\mathcal P_s\subseteq\mathcal P$ be nonempty.
Suppose \maincref{ass:main-process} holds on $\mathcal P_s$,
with nonempty binding classes satisfying
\mainequationref{eq:uniform-secants}.
Then
\[
 \sup_{P\in\mathcal P_s}
 d_K(\mathcal L_P(T_n(M_P)),\mathcal L(T_P^0))\to0.
\]
If random direction classes $A_{n,P}$ associated with selected pair sets satisfy
$d_{H,\rho_P}(\overline{A_{n,P}}^{\,L^2(P)},\mathcal H_{P,0})
=o_{\mathcal P_s}(1)$, let $S_n^\#$ be their clipped multiplier
supremum, dividing each comparison by its empirical standard deviation,
and let $c_n=q_{1-\alpha}(S_n^\#\mid Y_1,\ldots,Y_n)$. Then
\[
 \begin{aligned}
 d_K(\mathcal L_\xi(S_n^\#\mid Y_1,\ldots,Y_n),
       \mathcal L((T_P^0)_+))&=o_{\mathcal P_s}(1),\\
 c_n-c_{P,\alpha}^0&=o_{\mathcal P_s}(1),\\
 \sup_{P\in\mathcal P_s}
 \left|\pr_P(M_P\subseteq\mathcal C_n(c_n))-(1-\alpha)\right|&\to0.
 \end{aligned}
\]
\end{lemma}
\begin{proof}
By \eqref{eq:small-slack-risk} and
$\mathcal H_{P,0}\subseteq\mathcal H_{P,\eta}$,
\begin{equation}
 d_{H,\rho_P}(\mathcal H_{P,\eta},\mathcal H_{P,0})
 \le\omega_P(4\eta(1+\eta)F_P^\star).
 \label{eq:uniform-binding-localization}
\end{equation}
By \mainequationref{eq:uniform-secants}, the right side tends to zero
uniformly over $\mathcal P_s$.
Equations \eqrefs{eq:binding-lower,eq:binding-upper}, uniform
equicontinuity and studentization, and
\eqrefs{eq:uniform-donsker,eq:uniform-anticoncentration}
prove the first assertion.
For the selected classes, Hausdorff convergence and conditional
equicontinuity compare their suprema with finite nets of
$\mathcal H_{P,0}$ as in \eqref{eq:empirical-modulus}.
The studentization error is bounded by
\[
 \|\mathbb G_n^\#\|_{\mathcal H_P}
 \sup_{\mathcal I_P}
     \left|\frac{\sigma_P}{\wh\sigma_n}-1\right|
 =o_{\mathcal P_s,\xi}(1).
\]
The uniform multiplier approximation and \eqref{eq:uniform-anticoncentration}
give conditional Kolmogorov convergence before clipping;
$x\mapsto x_+$ cannot increase this distance.
Equation \eqref{eq:uniform-quantile-separation} gives critical-value
consistency. On the common event of
\cref{lem:zero-variance,lem:relative-variance}, inversion is equivalent to
$T_n(M_P)\le c_n$. On $|c_n-c_{P,\alpha}^0|\le\epsilon$, the first assertion
gives, for every $\epsilon>0$,
\[
 \pr(T_P^0\le c_{P,\alpha}^0-\epsilon)+o(1)
 \le \pr_P(M_P\subseteq\mathcal C_n(c_n))\le \pr(T_P^0\le c_{P,\alpha}^0+\epsilon)+o(1),
\]
uniformly over $\mathcal P_s$. Let $\epsilon\downarrow0$.
\end{proof}
\subsection{Global calibration}
\begin{proof}[\proofheading{\mainCref{thm:calibration}}]
Apply \cref{lem:uniform-quantiles} with
$\mathcal D_{n,P}=\mathcal I_P$. On the defining event of
\cref{lem:uniform-quantiles}, its clipped multiplier supremum is exactly $S_n^{\#,G}$.
Symmetry and Gaussian path continuity give
$T_P^G=\sup_{\overline{\mathcal H_P}^{\,L^2(P)}}\mathbb G_P
       =\sup_{\overline{\mathcal H_P}^{\,L^2(P)}}|\mathbb G_P|$.
The conditional approximation and \eqref{eq:uniform-quantile-separation} give
$c_n^G-c_{P,\alpha}^G=o_{\mathcal P}(1)$.

Since $\lambda_P(m,r)\ge0$,
\mainequationref{eq:studentized-mean-decomposition} gives
$(T_n(M_P))_+\le(U_{n,P})_+$. Zero-variance comparisons satisfy
inversion simultaneously by \cref{lem:zero-variance}. Thus, with
probability tending to one uniformly,
\begin{equation}
 \{M_P\subseteq\mathcal C_n(c)\}
       =\{T_n(M_P)\le c\},\qquad c\ge0,
 \label{eq:coverage-event-on-good-set}
\end{equation}
including data-dependent $c$. The final assertion of
\cref{lem:uniform-quantiles} proves the uniform coverage bound.
\end{proof}
\begin{lemma}[Relative relaxation and localization]
\label{lem:bounded-ball-localization}
Under \maincref{ass:main-process} and
\mainequationref{eq:adaptive-tuning-main}, put
$\widehat\delta_n=(2b_n+12\tau_n)F_n^\star/n^{1/2}$.
Uniformly over $P\in\mathcal P$, with probability tending to one,
\begin{equation}
 M_P\subseteq\widehat M_n
       \subseteq M_P(2b_nF_n^\star/n^{1/2}),\qquad
 \overline{\mathcal H_{P,n}^{\mathrm L}}^{\,L^2(P)}
       \subseteq\mathcal S_P(\widehat\delta_n).
 \label{eq:relaxed-tube-inclusions}
\end{equation}
Moreover, $\widehat\delta_n/F_P^\star=o_{\mathcal P}(1)$.
For each $P$, $F_n^\star\to_pF_P^\star>0$ and
$M_P(b_nF_P^\star/(4n^{1/2}))\subseteq\widehat M_n$ with probability
tending to one.
\end{lemma}
\begin{proof}
Put $k_n=\|\mathbb G_n\|_{\mathcal H_P}/n^{1/2}$, $r_n=b_n/n^{1/2}$ and
$t_n=\tau_n/n^{1/2}$ within this proof. Since $F_P^\star=0$ would
make $P$ a point mass and $\mathcal H_P$ empty,
$F_P^\star>0$ under \maincref{ass:main-process}. Work on the common
support and relative-error event of \cref{lem:relative-variance}, where
\begin{equation}
 |F_n(q)-F_n(r)-\{F_P(q)-F_P(r)\}|
 \le2k_n\wh\sigma_n(q,r)\qquad(q,r\in\cQ).
 \label{eq:sample-centered-error}
\end{equation}
Zero-variance pairs are included by \cref{lem:zero-variance}.
Since $\|\mathbb G_n\|_{\mathcal H_P}=O_{\mathcal P}(1)$ and
\mainequationref{eq:adaptive-tuning-main} holds, all restrictions
on $k_n$ below hold with probability tending to one uniformly.

For $m\in M_P$, compare with $r_j$ such that
$F_n(r_j)\downarrow F_n^\star$. Equations
\eqrefs{eq:sample-centered-error,eq:loss-risk-scale}, with $R=P_n$, give
\[
 F_n(m)-F_n(r_j)\le4k_n\{F_n(m)+F_n(r_j)\}.
\]
Let $j\to\infty$ and rearrange:
\begin{equation}
 \sup_{m\in M_P}\{F_n(m)-F_n^\star\}
 \le\frac{8k_n}{1-4k_n}F_n^\star\le r_nF_n^\star
 \quad\left(k_n\le\frac{r_n}{8+4r_n}\right).
 \label{eq:self-normalized-anchor}
\end{equation}

On the retention event, $q\in\widehat M_n$ and $m\in M_P$ have
$F_n(q),F_n(m)\le(1+r_n)F_n^\star$. Hence
\[
 F_P(q)-F_P^\star\le\{r_n+8k_n(1+r_n)\}F_n^\star
 \le2r_nF_n^\star
 \quad\left(k_n\le\frac{r_n}{8(1+r_n)}\right).
\]

For $(q,r)\in\widehat{\mathcal B}_n^{\mathrm L}$, the selection
inequality and \eqref{eq:loss-risk-scale} imply
\[
 F_n(r)\le\frac{1+2t_n}{1-2t_n}F_n(q),\qquad
 \wh\sigma_n(q,r)\le6F_n^\star
 \quad(r_n,t_n\le1/8).
\]
Consequently, on $2k_n\le t_n$,
\[
 F_P(r)-F_P(q)\le(t_n+2k_n)\wh\sigma_n(q,r)
 \le12t_nF_n^\star.
\]
Both endpoints lie in $M_P(\widehat\delta_n)$; taking the direction
closure gives \eqref{eq:relaxed-tube-inclusions}.

Choose a fixed $m_P\in M_P$. Since $F_n^\star\le F_n(m_P)$,
\[
 \sup_P\pr_P(F_n^\star/F_P^\star>L)
 \le\sup_P\frac{\E_PF_n(m_P)}{LF_P^\star}=L^{-1}.
\]
Thus $\widehat\delta_n/F_P^\star=o_{\mathcal P}(1)$.
For fixed $P$, \eqref{eq:loss-risk-scale}, now with $R=P$, gives
\[
 \begin{aligned}
 F_n(q)-F_n(m_P)&\ge(1-2k_n)\{F_P(q)-F_P^\star\}-4k_nF_P^\star,\\
 F_n(m_P)-4k_nF_P^\star&\le F_n^\star\le F_n(m_P),\\
 F_n(q)-F_n^\star&\le(1+2k_n)\{F_P(q)-F_P^\star\}+8k_nF_P^\star.
 \end{aligned}
\]
For $2k_n<1$, the second line and the law of large numbers give
$F_n^\star\to_pF_P^\star$. Since $k_n/r_n\to_p0$, the last line
is at most $r_nF_P^\star/2$ on $M_P(r_nF_P^\star/4)$ with probability
tending to one. Together with $F_n^\star\ge3F_P^\star/4$, this
proves the inner inclusion.
\end{proof}
\subsection{Localized calibration}
\label{sec:supp-zero-cutoff}
\begin{proof}[\proofheading{\mainCref{thm:zero-cutoff-main}}]
By \cref{lem:zero-variance,lem:relative-variance,lem:bounded-ball-localization},
with probability tending to one uniformly,
\begin{equation}
 \begin{gathered}
 M_P\subseteq\widehat M_n,\qquad \|\mathbb G_n\|_{\mathcal H_P}\le\tau_n/32,\\
 \{(q,r):\wh\sigma_n(q,r)>0\}=\mathcal I_P,\qquad
 \frac12\le\frac{\sigma_P(q,r)}{\wh\sigma_n(q,r)}\le2
 \quad((q,r)\in\mathcal I_P).
 \end{gathered}
 \label{eq:proof-gms-event}
\end{equation}
Use \cref{lem:uniform-quantiles} with the deterministic pair set
\[
 \mathcal D_{n,P}=
 \{(m,r):m\in M_P,\ (m,r)\in\mathcal I_P,\,
                  \lambda_P(m,r)\le\tau_n/(16n^{1/2})\}.
\]
On \eqref{eq:proof-gms-event}, each pair in $\mathcal D_{n,P}$ satisfies
\[
 T_n(m,r)\ge-2\{\|\mathbb G_n\|_{\mathcal H_P}+\tau_n/16\}
       \ge-3\tau_n/16>-\tau_n.
\]
Consequently,
\begin{equation}
 \mathcal A_P(\tau_n/(16n^{1/2}))\subseteq\mathcal H_{P,n}^{\mathrm L},
 \qquad
 \mathcal H_{P,0}\subseteq
       \overline{\mathcal H_{P,n}^{\mathrm L}}^{\,L^2(P)}.
 \label{eq:proof-binding-retention}
\end{equation}
Outside $\mathcal D_{n,P}$,
$\mathbb G_nh_{m,r}-n^{1/2}\lambda_P(m,r)<0$.
If $\mathcal D_{n,P}\ne\varnothing$, the statistics in
\cref{lem:uniform-quantiles} therefore satisfy
\[
 (T_n(M_P))_+\le(U_{n,P})_+,\qquad
 S_n^{\#,\mathrm L}\ge(V_{n,P}^\#)_+.
\]
Quantile monotonicity, \cref{lem:uniform-quantiles} and
\eqref{eq:coverage-event-on-good-set} give coverage.
If $\mathcal D_{n,P}=\varnothing$, every mean-anchored statistic is
negative on \eqref{eq:proof-gms-event}, so the nonnegative cutoff
covers $M_P$. Both cases hold along arbitrary sequences of laws in $\mathcal P$.

For nonempty binding classes,
\eqrefs{eq:proof-binding-retention,eq:relaxed-tube-inclusions} imply,
with probability tending to one uniformly,
\begin{equation}
 d_{H,\rho_P}(\overline{\mathcal H_{P,n}^{\mathrm L}}^{\,L^2(P)},\mathcal H_{P,0})
 \le\omega_P(\widehat\delta_n).
 \label{eq:selected-binding-distance}
\end{equation}
Since $\widehat\delta_n/F_P^\star=o_{\mathcal P}(1)$,
\mainequationref{eq:uniform-secants} and
\cref{lem:uniform-binding-transfer} give both uniform exactness conclusions.
\end{proof}

\subsection{Intrinsic circle}
\label{sec:supp-circle-theory}
\begin{proof}[\proofheading{\mainCref{cor:circle-exactness-main}}]
Write $\operatorname{wrap}(x)$ for the representative of $x$ in
$[-\pi,\pi)$. We verify \maincref{ass:main-process}, nonemptiness of
$\mathcal H_{P,0}$ and \mainequationref{eq:uniform-secants}, then apply
\maincref{thm:zero-cutoff-main}.

\noindent\emph{Empirical-process regularity.}
For $q\ne r$, let
$U_{q,r}=\{d^2(\cdot,q)-d^2(\cdot,r)\}/d(q,r)$.
If $P$ has exactly two support points, contrasting those points shows
that $\mathcal H_P$ is nonempty. Since $L_0^2(P)$ is one-dimensional,
its normalized directions are two opposite functions. This finite
class has a bounded envelope on the support and satisfies
\maincref{ass:main-process} by the scalar CLT.

Suppose that the support contains at least three points. After rotation
and, if necessary, interchanging the endpoints, take $q=0$ and
$r=t\in(0,\pi]$. For $x\in[-\pi,\pi)$,
\begin{equation}
 d^2(x,0)-d^2(x,t)=
 \begin{cases}
 -2(2\pi-t)x-(2\pi-t)^2,&-\pi\le x<-\pi+t,\\
 2tx-t^2,&-\pi+t\le x<\pi.
 \end{cases}
 \label{eq:circle-loss-contrast}
\end{equation}
Thus $|U_{q,r}|\le2\pi$ and every level set of $U_{q,r}$ contains
at most two points. In fact,
\begin{equation}
 \inf_{q\ne r}\Var_P(U_{q,r})>0.
 \label{eq:circle-positive-rescaled-variance}
\end{equation}
Otherwise, choose pairs with variances tending to zero and pass to
$q_j\to q$ and $r_j\to r$. If $q\ne r$, the rescaled contrasts
converge uniformly to $U_{q,r}$. If $q=r$, pass to a fixed orientation.
Away from the antipode of $q$, they converge pointwise to a signed
version of $2\operatorname{wrap}(y-q)$; a further subsequence makes
their bounded values at the antipode converge as well. Each limit
has level sets of at most two points. Bounded convergence gives
$L^2(P)$ convergence, and zero limiting variance would concentrate
$P$ on at most two points, a contradiction.

By \eqref{eq:circle-positive-rescaled-variance},
$h_{q,r}=(U_{q,r}-PU_{q,r})/\{\Var_P(U_{q,r})\}^{1/2}$ has a bounded
envelope. In a fixed angular chart each normalized contrast is affine
on at most three intervals, separated by the two antipodal cuts.
Its subgraph is a union of at most three intersections of vertical
strips and half-planes. These classes have a common polynomial
VC-subgraph entropy bound
\citep[Lemma~2.6.19 and Theorem~2.6.7]{VaartWellner2023}.
Pointwise separability follows from continuity in $(q,r)$ on the
positive-variance index set. Apply \cref{lem:uniform-entropy} to
$\{P\}$ to obtain \maincref{ass:main-process}.

\medskip\noindent\emph{Fixed-law exactness.}
Compactness gives a nonempty compact mean set and $F_P^\star>0$.
At any $m\in M_P$, bounded one-sided derivatives of squared distance
give
\[
 F'_{P,+}(m)=F'_{P,-}(m)-4\pi P(Y=m+\pi).
\]
A local minimum requires $F'_{P,-}(m)\le0\le F'_{P,+}(m)$, so the
antipode has no atom and both derivatives vanish. Consequently,
\[
 \psi_m(Y)=2\operatorname{wrap}(Y-m),\qquad
 P\psi_m=0,\qquad \|\psi_m\|_2^2=4F_P^\star>0.
\]
By \eqref{eq:circle-loss-contrast}, for $|u|,|v|\le\delta<\pi/2$,
the uncentred contrast equals
$(v-u)\psi_m(Y)+u^2-v^2$ outside the radius-$\delta$ arc about the
antipode. The wrapping correction has absolute value at most
$4\pi|v-u|\mathds1\{d(Y,m+\pi)\le\delta\}$.
Centring is an $L^2(P)$ contraction, hence
\begin{equation}
 \sup_{\substack{|u|,|v|\le\delta\\u\ne v}}
 \frac{\|\ell_{m+u,m+v}-(v-u)\psi_m\|_2}{|v-u|}
 \le4\pi P\{d(Y,m+\pi)\le\delta\}^{1/2}\longrightarrow0.
 \label{eq:circle-strict-secant}
\end{equation}
Normalizing the contrasts multiplies this error bound by at most
$2/\|\psi_m\|_2$. Since $F_P(m+t)-F_P^\star=o(|t|)$ and
$\sigma_P(m,m+t)=|t|\|\psi_m\|_2+o(|t|)$, comparisons anchored at
$m$ with $t\to0$ from either side have vanishing standardized slack
and limiting directions $\pm\psi_m/\|\psi_m\|_2$. Both belong to
$\mathcal H_{P,0}$, which is therefore nonempty.

For two-point support these are all the normalized directions, so
$\omega_P(\delta)=0$ for every $\delta>0$. Otherwise, if
$\omega_P(\delta)\not\to0$, approximation by actual pairs in the
definition of $\mathcal S_P(\delta)$ gives $\delta_j\downarrow0$ and
positive-variance pairs $q_j,r_j\in M_P(\delta_j)$ whose normalized
contrasts remain a fixed positive distance from the closed binding
class. A subsequence has limits $m,m'\in M_P$.
If $m\ne m'$, \eqref{eq:circle-loss-contrast} gives
$\sigma_P(m,m')>0$ and continuity yields
$h_{q_j,r_j}\to h_{m,m'}$ in $L^2(P)$. The limit is binding because
both endpoints minimize risk. If $m=m'$, \eqref{eq:circle-strict-secant}
puts the normalized contrasts arbitrarily close to one of the two
binding score directions at $m$. Both cases contradict the choice
of the pairs. Thus $\omega_P(\delta)\to0$. Since $F_P^\star>0$,
\mainequationref{eq:uniform-secants} holds for $\{P\}$, and
\maincref{thm:zero-cutoff-main} gives critical-value consistency and
exact whole-mean-set coverage.

\medskip\noindent\emph{Uniformity over the circular family.}
The density numerator in \mainequationref{eq:circle-family} is affine
in $\lambda$; its minima at $\lambda=-1/8$ and $\lambda=1$ are
$263/1536$ and $1/2$. The densities therefore have common positive
lower and finite upper bounds. For uniform circular probability
$\mu$, $P_\lambda\ge(263/1536)\mu$. Using
$\Var_P(U)=\inf_{a\in\mathbb R}P(U-a)^2$ and applying
\eqref{eq:circle-positive-rescaled-variance} to $\mu$ gives
\[
 \inf_{\lambda\in[-1/8,1]}\inf_{q\ne r}\Var_{P_\lambda}(U_{q,r})
 \ge\frac{263}{1536}\inf_{q\ne r}\Var_\mu(U_{q,r})>0.
\]
The common bounded envelope and the same entropy bound verify
\maincref{ass:main-process} uniformly by \cref{lem:uniform-entropy}.
The density bounds also give
$\inf_\lambda F_{P_\lambda}^\star>0$ and make the error in
\eqref{eq:circle-strict-secant} of order $O(\delta^{1/2})$, uniformly
in $\lambda$ and $m\in M_{P_\lambda}$.

Direct integration of \mainequationref{eq:circle-family} gives
\begin{equation}
 F_{P_\lambda}(t)=\pi^2/3-1+\lambda(1-\cos t)
                         +\tfrac23(1-\lambda)(1-\cos t)^2.
 \label{eq:circle-family-risk}
\end{equation}
Minimizing this quadratic in $1-\cos t$ yields
\[
 M_{P_\lambda}=\begin{cases}
 \{-t_\lambda,t_\lambda\},&\lambda<0,\\
 \{0\},&\lambda\ge0,
 \end{cases}\qquad
 t_\lambda=\arccos\!\left(\frac{4-\lambda}{4(1-\lambda)}\right)
 \quad(\lambda<0).
\]
These sets are Hausdorff-continuous in $\lambda$, including at zero.
The risks and their minima are continuous in the parameter.
If \mainequationref{eq:uniform-secants} failed, approximation by
actual pairs would give $u_j\downarrow0$, parameters $\lambda_j$ and
positive-variance pairs
$q_j,r_j\in M_{P_{\lambda_j}}(u_jF_{P_{\lambda_j}}^\star)$
whose normalized contrasts remain a fixed positive distance from
$\mathcal H_{P_{\lambda_j},0}$. Pass to a subsequence with
$\lambda_j\to\lambda_*$, $q_j\to m$ and $r_j\to m'$.
Joint continuity and $F_{P_\lambda}^\star\le\pi^2$ imply
$m,m'\in M_{P_{\lambda_*}}$.
If $m\ne m'$, Hausdorff continuity provides means $m_j,m'_j$ for
$P_{\lambda_j}$ converging to $m,m'$. The limiting contrast has
positive variance, and continuity of the bounded losses and densities
gives
\[
 \|h_{q_j,r_j}^{(P_{\lambda_j})}
       -h_{m_j,m'_j}^{(P_{\lambda_j})}\|_{P_{\lambda_j},2}\to0,
\]
where both directions are centred and normalized under
$P_{\lambda_j}$. The second is binding, a contradiction.
If $m=m'$, choose $m_j\in M_{P_{\lambda_j}}$ with $m_j\to m$.
Both endpoints lie in a shrinking chart about $m_j$.
The uniform bound in \eqref{eq:circle-strict-secant} and
$\|\psi_{m_j}\|_{P_{\lambda_j},2}^2=4F_{P_{\lambda_j}}^\star$
bounded away from zero give
\[
 \min_{\varepsilon\in\{-1,1\}}
 \left\|h_{q_j,r_j}^{(P_{\lambda_j})}
   -\varepsilon\frac{\psi_{m_j}}{\|\psi_{m_j}\|_{P_{\lambda_j},2}}
 \right\|_{P_{\lambda_j},2}\to0.
\]
Both score directions are binding by the fixed-law argument, again
a contradiction. This proves \mainequationref{eq:uniform-secants}
for the whole family. Apply \maincref{thm:zero-cutoff-main} to obtain
critical-value consistency and exact whole-mean-set coverage,
uniformly over the circular family.
\end{proof}

\section{Sufficient conditions}
\subsection{Empirical-process conditions}
\label{sec:supp-entropy}
\begin{lemma}[Uniform entropy criterion]
\label{lem:uniform-entropy}
Suppose the nonempty classes $\mathcal H_P$ have measurable envelopes
$H_P$ satisfying the uniform square-integrability
condition in \maincref{ass:main-process}. Suppose also
\begin{equation}
 \int_0^1\sup_{P\in\mathcal P}\sup_\nu
 \left\{1+\log N\bigl(
 u\|H_P\|_{\nu,2},\mathcal H_P,L^2(\nu)\bigr)\right\}^{1/2}\,du
 <\infty,
 \label{eq:uniform-entropy-criterion}
\end{equation}
where $\nu$ ranges over finitely supported probability measures on
a $P$-full set on which the envelope and separability properties
hold, with finite positive envelope norm. Then \maincref{ass:main-process} holds.
A common polynomial bound $(A/u)^v$ on these covering numbers suffices.
\end{lemma}
\begin{proof}
The entropy bounds in \citet[Theorem~2.3]{SheehyWellner1992}
and \citet[Theorem~2.8.3]{VaartWellner2023} apply to each row
$\mathcal H_P$ with constants controlled by
\eqref{eq:uniform-entropy-criterion}.
For the truncated squared-increment class, symmetrization and
chaining give, for each $L>0$,
\[
 \sup_P\E_P\left[\sup_{h,g\in\mathcal H_P}
 \left|\frac1n\sum_i(h-g)^2(Y_i)-P(h-g)^2\right|\right]\le\frac{C_L}{n^{1/2}}
       +8\sup_P \E_P[H_P(Y)^2\mathds1\{H_P(Y)>L\}],
\]
where $C_L<\infty$ depends only on $L$ and the common entropy integral.
Let $n\to\infty$, then $L\to\infty$. The empirical squared-increment
metric approaches $\rho_P^2$ uniformly; local chaining and
\eqref{eq:uniform-entropy-criterion} yield
\[
 \lim_{\delta\downarrow0}\limsup_n\sup_P
 \E_P\left[\sup_{\substack{h,g\in\mathcal H_P\\\rho_P(h,g)<\delta}}
                 |\mathbb G_n(h-g)|\right]=0.
\]
Gaussian chaining gives \eqref{eq:uniform-pregaussian} and finite
$\rho_P$-nets of uniformly bounded cardinality. On such a net,
\[
 \sup_P\max_{j\le N}
 \E_P[h_j(Y)^2\mathds1\{|h_j(Y)|>\epsilon n^{1/2}\}]
 \le\sup_P \E_P[H_P(Y)^2\mathds1\{H_P(Y)>\epsilon n^{1/2}\}]
 \longrightarrow0.
\]
This uniform Lindeberg condition and the two modulus bounds give
\eqref{eq:uniform-donsker} by the finite-net argument of
\citet[Theorem~2.8.2]{VaartWellner2023}.
\end{proof}

Let $\mathcal V_P\coloneqq
\operatorname{span}(\{\ell_{q,r}:(q,r)\in\mathcal I_P\})
\subseteq L_0^2(P)$ and $R_P\coloneqq\dim(\mathcal V_P)$.
\begin{lemma}[Finite-rank and H\"older criteria]
\label{lem:ep-sufficient}
Fix $P$. Either condition below implies \maincref{ass:main-process} with
$\mathcal P=\{P\}$.
\begin{enumerate}
\item The centred-loss rank satisfies $1\leq R_P<\infty$.
\item There are a nonempty compact $\Theta\subset\mathbb R^k$,
$0<\gamma\leq1$, and a representation
$\mathcal H_P=\{h_\theta:\theta\in\Theta\}$ such that, on a common
full-probability set, for some nonnegative $L\in L^2(P)$,
$|h_\theta-h_{\theta'}|\leq L\|\theta-\theta'\|^\gamma$.
\end{enumerate}
\end{lemma}
\begin{proof}
Under (i), choose an orthonormal basis $\psi_P$ of $\mathcal V_P$.
For every $h\in\mathcal H_P$,
\[
 \begin{aligned}
 &h=a^{\mathrm T}\psi_P,\quad\|a\|=1,\quad
 P\psi_P=0,\quad P\psi_P\psi_P^{\mathrm T}=I_{R_P},\\
 &\mathbb G_nh=a^{\mathrm T}\mathbb G_n\psi_P,\quad
 |h|\le\|\psi_P\|,\quad P\|\psi_P\|^2=R_P<\infty.
 \end{aligned}
\]
The multivariate CLT and the 1-Lipschitz linear map
to the coefficient-indexed process give Donsker convergence.
Choose the basis as finite linear combinations of continuous contrasts;
all identities then hold on the common support, and a countable dense
coefficient set gives separability.

Under (ii), a $\delta$-net $\theta_1,\ldots,\theta_N$ gives brackets
$[h_{\theta_j}-\delta^\gamma L,\,
h_{\theta_j}+\delta^\gamma L]$ of $L^2(P)$ width
$2\delta^\gamma\|L\|_2$.
For $\|L\|_2>0$, compactness in $\mathbb R^k$ implies
\[
 N_{[]}(\epsilon,\mathcal H_P,L^2(P))
 \le C\left(1+\frac{\|L\|_2}{\epsilon}\right)^{k/\gamma},
 \qquad
 \int_0^1\left\{\log N_{[]}(\epsilon,\mathcal H_P,L^2(P))\right\}^{1/2}\,d\epsilon
 <\infty.
\]
If $\|L\|_2=0$, the class is a singleton in $L^2(P)$.
For $\theta_0\in\Theta$, the envelope $|h_{\theta_0}|+L\operatorname{diam}(\Theta)^\gamma$
is square-integrable. Apply the bracketing Donsker theorem
\citep[Theorem~2.5.6]{VaartWellner2023}. The common pointwise H\"older
bound gives a countable separable version.
\end{proof}

The finite-rank criterion is uniform if $1\le R_P\le R<\infty$ and
the orthonormal-basis envelopes $\|\psi_P\|$ are square integrable
uniformly in $P$. Indeed, the coefficient balls have a common
polynomial covering bound under every $L^2(\nu)$, using the envelope
$\|\psi_P\|$. The H\"older criterion is uniform if $\Theta$, its
dimension and $\gamma$ can be chosen in common, and, for some fixed $\theta_0\in\Theta$,
$|h_{P,\theta_0}|+\{1+\operatorname{diam}(\Theta)^\gamma\}L_P$
is a uniformly square-integrable envelope. A Euclidean parameter net
then gives a common polynomial covering bound. Apply
\cref{lem:uniform-entropy} in either case.

If $\operatorname{supp}(P)=\{y_1,\ldots,y_s\}$ with positive masses,
then $L_0^2(P)$ has dimension $s-1$ and hence
$R_P\leq s-1$, regardless of whether $\cQ$ is finite.
When $s\geq2$, contrasting two distinct support points shows
$R_P\geq1$, so the process condition holds. A single support point
has no normalized directions and is covered by the zero-variance case.
If every $P$ has finite support of at least two points and
every support mass is at least $p_*>0$, then $|\operatorname{supp}(P)|\le p_*^{-1}$ and
$Ph^2=1$ gives $|h(y)|\le p_*^{-1/2}$ at every support point. Thus these models
satisfy \maincref{ass:main-process} uniformly.

\subsection{Sufficient conditions for direction recovery}
\label{sec:supp-directions}
On a finite metric space, the binding class consists exactly of the
normalized positive-variance contrasts between population minimizers;
in particular, it is empty when the population mean is unique.

\begin{proposition}[Finitely many stationary means]
\label{prop:finite-stationary-means}
Suppose $M_P=\{m_1,\ldots,m_J\}$ with $J\ge1$ and
$F_P(q_l)-F_P^\star\to0$ implies $d(q_l,M_P)\to0$ for every sequence $(q_l)$.
For each $j$, let $\phi_j$ be a homeomorphism from a neighbourhood of
$0\in\mathbb R^{d_j}$, $d_j\ge1$, onto a neighbourhood of $m_j$, with
$\phi_j(0)=m_j$. Suppose $F_P\circ\phi_j$ is differentiable at zero
with derivative zero and, for an injective linear map
$A_j:\mathbb R^{d_j}\to L_0^2(P)$,
\begin{equation}
 \sup_{\substack{u,v\in B(0,\delta)\\u\ne v}}
 \frac{\|\Phi_P(\phi_j(u))-\Phi_P(\phi_j(v))-A_j(v-u)\|_2}
 {\|u-v\|}\longrightarrow0.
 \label{eq:strict-secant-differentiability}
\end{equation}
Assume $\sigma_P(m_j,m_k)>0$ whenever $j\ne k$. Then
\begin{equation}
 \mathcal H_{P,0}=
 \bigcup_{j=1}^J\left\{\frac{A_jw}{\|A_jw\|_2}:\|w\|=1\right\}
 \cup\{h_{m_j,m_k}:j\ne k\},
 \label{eq:finite-stationary-binding}
\end{equation}
and $\omega_P(\delta)\to0$ as $\delta\downarrow0$.
\end{proposition}
\begin{proof}
Injectivity in finite dimension gives
$\inf_{\|w\|=1}\|A_jw\|_2>0$; the right-hand side of
\eqref{eq:finite-stationary-binding} is compact.
For any sequence of positive-variance pairs with vanishing excess
risks, isolation and finiteness give a subsequence whose endpoints
converge to $m_j,m_k$.
If $j=k$, the strict secant expansion puts the normalized contrast at
vanishing distance from the $j$th unit-direction class.
If $j\ne k$, continuity of $\Phi_P$ and the positive cross-mean variance
give $h_{q,r}\to h_{m_j,m_k}$.
Hence all low-risk directions approach the compact set on the
right-hand side of \eqref{eq:finite-stationary-binding}.

For a binding direction choose mean-anchored comparisons with slack
tending to zero and directions tending to it. Equation
\eqref{eq:small-slack-risk} makes their second endpoints near-minimizers;
finiteness fixes their first endpoint along a subsequence.
The same two cases give inclusion in this compact set.
Conversely, cross-mean comparisons have zero slack, and for each unit $w$,
\[
 F_P(\phi_j(tw))-F_P^\star=o(t),\qquad
 \ell_{m_j,\phi_j(tw)}=tA_jw+o_{L^2(P)}(t),\qquad
 \sigma_P(m_j,\phi_j(tw))=t\|A_jw\|_2+o(t).
\]
Their slack tends to zero and their directions converge to
$A_jw/\|A_jw\|_2$. The two inclusions prove
\eqref{eq:finite-stationary-binding}; the low-risk approximation gives
$\omega_P(\delta)\to0$, also after taking closures.
\end{proof}
Continuous Fr\'echet differentiability of $\Phi_P\circ\phi_j$ implies
\eqref{eq:strict-secant-differentiability} with
$A_j=-D(\Phi_P\circ\phi_j)(0)$, since the increment formula gives
\[
 \frac{\|\Phi_P(\phi_j(u))-\Phi_P(\phi_j(v))-A_j(v-u)\|_2}{\|u-v\|}
 \le\sup_{\|x\|\le\delta}
 \|D(\Phi_P\circ\phi_j)(x)+A_j\|_{\mathrm{op}}\to0.
\]

\section{Wald recovery and confidence radii}
\subsection{Euclidean means}
\label{sec:supp-euclidean}
\begin{proof}[\proofheading{\mainCref{prop:euclidean}}]
For $a\in\ker(\Sigma)$,
$\E_P[\{a^{\mathrm T}(Y-\mu)\}^2]=a^{\mathrm T}\Sigma a=0$.
Applying this to a finite basis gives
$Y-\mu\in\operatorname{Range}(\Sigma)$ a.s..
The centred contrast identity in \maincref{sec:geometry} gives
$\sigma_P^2(q,s)=4(s-q)^{\mathrm T}\Sigma(s-q)$.
For $v^{\mathrm T}\Sigma v>0$,
\[
 h_{\mu,\mu+tv}
 =\frac{v^{\mathrm T}(Y-\mu)}{(v^{\mathrm T}\Sigma v)^{1/2}},
 \qquad
 \lambda_P(\mu,\mu+tv)
 =\frac{t\|v\|^2}{2(v^{\mathrm T}\Sigma v)^{1/2}}\to0.
\]
These generate $\mathcal H_P$, proving the binding-class equality
and the $\chi_r$ supremum law.

For $A\succeq0$ of positive rank,
\begin{equation}
 \sup_{v:v^{\mathrm T}Av>0}
 \frac{-v^{\mathrm T}x}{(v^{\mathrm T}Av)^{1/2}}
 =
 \begin{cases}
 (x^{\mathrm T}A^+x)^{1/2},&x\in\operatorname{Range}(A),\\
 +\infty,&x\notin\operatorname{Range}(A).
 \end{cases}
 \label{eq:quadratic-projection-identity}
\end{equation}
If $x\in\operatorname{Range}(A)$, Cauchy--Schwarz gives the bound,
with equality at $v=-A^+x$ when $x\ne0$. For $x=0$ it is zero.
Otherwise write $x=x_R+x_0$, where $x_R\in\operatorname{Range}(A)$
and $0\ne x_0\in\ker(A)$. Taking $v=u-a x_0$ with
$u^{\mathrm T}Au>0$ and $a\to\infty$ makes the numerator unbounded
without changing the denominator.

For $s=q+tv$, $t>0$, direct expansion gives
\[
 F_n(q)-F_n(q+tv)=-2t\,v^{\mathrm T}(q-\bar Y_n)-t^2\|v\|^2,\qquad
 \wh\sigma_n(q,q+tv)=2t(v^{\mathrm T}\wh\Sigma_n v)^{1/2}.
\]
Let $t\downarrow0$ and use \eqref{eq:quadratic-projection-identity}
with $A=\wh\Sigma_n$ and $x=q-\bar Y_n$. This proves
\mainequationref{eq:euclidean-statistic}.
If $q-\bar Y_n\in\operatorname{Range}(\wh\Sigma_n)$ and
$\wh\Sigma_n v=0$, the zero-variance numerator is
$-t^2\|v\|^2\le0$. Hence such comparisons impose no additional
restriction on that affine support.

Conditionally on the observations,
\[
 \frac1{n^{1/2}}\sum_i\xi_i(Y_i-\bar Y_n)\sim N_d(0,\wh\Sigma_n).
\]
Taking all standardized projections of this vector and applying
\eqref{eq:quadratic-projection-identity} gives
\[
 \mathcal L_\xi(S_n^{\#,G}\mid Y_1,\ldots,Y_n)=\chi_{\wh r_n},
 \qquad c_n^G=(\chi^2_{\wh r_n,1-\alpha})^{1/2}.
\]
If $\wh r_n=0$, all observations coincide,
both cutoffs are zero, and inversion gives
$\{\bar Y_n\}$ directly.

For localization, $\bar Y_n\in\widehat M_n$ for every nonnegative
relaxation. For $v^{\mathrm T}\wh\Sigma_n v>0$,
\[
 T_n(\bar Y_n,\bar Y_n+tv)
 =-\frac{n^{1/2}\,t\|v\|^2}
          {2(v^{\mathrm T}\wh\Sigma_n v)^{1/2}}\to0,\qquad
 T_n^\#(\bar Y_n,\bar Y_n+tv)
 =\frac{v^{\mathrm T}n^{-1/2}\sum_i\xi_i(Y_i-\bar Y_n)}
          {(v^{\mathrm T}\wh\Sigma_n v)^{1/2}}.
\]
Since $\tau_n>0$, small $t$ retains every global multiplier direction. Thus
\[
 S_n^{\#,\mathrm L}=S_n^{\#,G},\qquad c_n^{\mathrm L}=c_n^G,\qquad C_n^{\mathrm L}=C_n^G.
\]

For uniform coverage, let $E_{n,P}$ be the defining event of
\cref{lem:uniform-quantiles} intersected with the probability-one event that all observations
lie in $\mu_P+\operatorname{Range}(\Sigma_P)$.
By \cref{lem:relative-variance}, $\sup_P\pr_P(E_{n,P}^c)\to0$.
The population and empirical variance formulas imply that, on $E_{n,P}$,
$\Sigma_P$ and $\wh\Sigma_n$ have the same kernel and range. Thus
$\wh r_n=r_P$ and $\bar Y_n-\mu_P\in\operatorname{Range}(\wh\Sigma_n)$,
proving uniform rank recovery.

Apply \cref{lem:uniform-quantiles}(i) with
$\mathcal D_{n,P}=\mathcal I_P$. Its Gaussian target is $\chi_{r_P}$, and
on $E_{n,P}$, \eqref{eq:quadratic-projection-identity} gives
\[
 U_{n,P}=\left\{n(\bar Y_n-\mu_P)^{\mathrm T}
 \wh\Sigma_n^+(\bar Y_n-\mu_P)\right\}^{1/2}.
\]
Put $\kappa_P=(\chi^2_{r_P,1-\alpha})^{1/2}$.
On $E_{n,P}$, both cutoffs equal $\kappa_P$, and
$\{\mu_P\in C_n^J\}=\{U_{n,P}\le\kappa_P\}$ for $J\in\{G,\mathrm L\}$.
Consequently,
\[
 \left|\pr_P(\mu_P\in C_n^J)-(1-\alpha)\right|
 \le\pr_P(E_{n,P}^c)+d_K(\mathcal L_P(U_{n,P}),\chi_{r_P})\to0
\]
uniformly over $P$ and $J$.

\end{proof}
\subsection{Local Wald recovery}
\label{sec:supp-wald-proof}
\label{sec:supp-wald-conditions}
Fix $P$ under the Polish-space and second-moment setup of
\maincref{sec:method}. No global Donsker condition is imposed.
Suppose $M_P=\{m\}$ and
\begin{equation}
 F_P(q_j)-F_P^\star\to0\quad\Longrightarrow\quad d(q_j,m)\to0.
 \label{eq:regular-risk-isolation}
\end{equation}
This isolation condition is automatic on a proper metric space:
the reverse triangle inequality makes risk sublevel sets bounded,
and compactness and uniqueness then give isolation.
Let $d\ge1$ and let $\phi:B(0,r_0)\subset\mathbb R^d\to\cQ$ be a homeomorphism onto
an open neighbourhood of $m$, with $\phi(0)=m$ and
$d(\phi(u),\phi(v))\le L\|u-v\|$.
Write $f_n=F_n\circ\phi$, $f=F_P\circ\phi$, and
$g(y,u)=d^2(y,\phi(u))$.
\begin{condition}[Local Wald regularity]
\label{cond:local-wald}
\leavevmode\par
\begin{enumerate}
\item[(W1)] For a measurable $\psi\in L^2(P;\mathbb R^d)$ with
$P\psi=0$ and $\Sigma=P(\psi\psi^{\mathrm T})>0$,
\begin{equation}
 \epsilon_A(\delta)\coloneqq
 \sup_{0<\|u\|\le\delta}
 \frac{\|\ell_{m,\phi(u)}-u^{\mathrm T}\psi\|_2}{\|u\|}
 \longrightarrow0\qquad(\delta\downarrow0).
 \label{eq:regular-anchored-expansion}
\end{equation}
\item[(W2)] For a symmetric positive-definite matrix $H$,
\begin{equation}
 f(u)-f(0)=\tfrac12u^{\mathrm T}Hu+o(\|u\|^2).
 \label{eq:regular-risk-expansion}
\end{equation}
\item[(W3)] For some $0<\delta_0<r_0$, the local comparison class
\[
 \mathcal H_{\mathrm{loc}}\coloneqq
 \{h_{\phi(u),\phi(v)}:u,v\in B(0,\delta_0),\
                         \sigma_P(\phi(u),\phi(v))>0\}
\]
is $P$-Donsker with a measurable square-integrable envelope.
With $\mathcal E\coloneqq
\{a^{\mathrm T}\psi/(a^{\mathrm T}\Sigma a)^{1/2}:\|a\|=1\}$,
assume also
\begin{equation}
 \sup_{\substack{u,v\in B(0,\delta)\\
                 \sigma_P(\phi(u),\phi(v))>0}}
 \inf_{e\in\mathcal E}\|h_{\phi(u),\phi(v)}-e\|_2
 \longrightarrow0\qquad(\delta\downarrow0).
 \label{eq:regular-direction-stability}
\end{equation}
\end{enumerate}
\end{condition}
\begin{remark}[Riemannian interpretation]
Suppose that $\cQ$ is a complete $d$-dimensional Riemannian manifold
with its geodesic distance and $\phi=\operatorname{Exp}_m$ is a normal
chart, using an orthonormal basis to identify $T_m\cQ$ with $\mathbb R^d$.
If $P(\operatorname{Cut}(m))=0$, first variation and dominated convergence
under the standing second-moment assumption give $P\psi=0$ and the
$L^2$ expansion in (W1) with $\psi(Y)=2\operatorname{Log}_m(Y)$.
Thus (W1) holds provided that
$\Sigma=4\Cov_P(\operatorname{Log}_m(Y))$ is positive definite.
If $F_P$ is twice differentiable at $m$ with positive-definite
Riemannian Hessian, (W2) holds with $H=\operatorname{Hess}F_P(m)$
in these coordinates.
These are the score and curvature quantities in regular intrinsic
Fr\'echet mean CLTs \citep{BhattacharyaPatrangenaru2005}; (W3) remains
an additional condition on normalized pairwise contrasts.
\end{remark}

Put $W_n=\mathbb G_n\psi$ and $V=H^{-1}\Sigma H^{-1}$.
When $M_n\ne\varnothing$, choose any measurable empirical mean
$\widehat m_n\in M_n$ and set
\[
 \widehat u_n\coloneqq
 \begin{cases}
 \phi^{-1}(\widehat m_n),&M_n\ne\varnothing,
                         \ \widehat m_n\in\phi(B(0,r_0)),\\
 0,&\text{otherwise}.
 \end{cases}
\]
For the Wald benchmark in \maincref{prop:wald-recovery}, let
$\widehat V_n=V+o_p(1)$ be symmetric positive definite.

\medskip\noindent\emph{Smooth sufficient conditions.}
Alternatively, suppose $g(Y,\cdot)$ is a.s.\ $C^2$
on a fixed chart ball. Write $s(y,u)=\nabla g(y,u)$,
$A(y,u)=\nabla^2g(y,u)$ and
$J_\delta(y)=\sup_{\|u\|\le\delta}\|A(y,u)-A(y,0)\|_{\mathrm{op}}$.
Assume
\begin{equation}
 P\|A(Y,0)\|_{\mathrm{op}}<\infty,\qquad
 PJ_\delta\to0,\qquad
 H=PA(Y,0)>0,\qquad \Sigma=\Var_P(s(Y,0))>0,
 \label{eq:regular-smooth-moments}
\end{equation}
and, for every $\eta>0$,
\begin{equation}
 \lim_{\delta\downarrow0}\limsup_n
 \pr_P\left(\E_\xi\left[1\wedge
 \sup_{\|u\|\le\delta}
 \|\mathbb G_n^\#(s(\cdot,u)-s(\cdot,0))\|\right]>\eta\right)=0.
 \label{eq:regular-score-multiplier}
\end{equation}
The chart bound gives
$\sup_u\|s(Y,u)\|\le2L\{d(Y,m)+Lr_0\}\eqqcolon G(Y)\in L^2(P)$.
Keeping \eqref{eq:regular-risk-isolation} and the chart assumptions,
\eqrefs{eq:regular-smooth-moments,eq:regular-score-multiplier} provide an
alternative to \cref{cond:local-wald}.

For example, \eqref{eq:regular-score-multiplier} follows if
$P(GJ)<\infty$, where
$J(Y)=\sup_{\|u\|\le\delta_1}\|A(Y,u)\|_{\mathrm{op}}$
for some fixed small $\delta_1$. Indeed,
$\|s(Y,u)-s(Y,v)\|\le\min(J(Y)\|u-v\|,2G(Y))$.
Conditional Gaussian increment distances are consequently bounded by
$\{2P_n(GJ)\|u-v\|\}^{1/2}$. Polynomial covering numbers of the
coordinate ball and the vanishing local $L^2(P_n)$ diameter give
\eqref{eq:regular-score-multiplier} by the Gaussian entropy integral.
On bounded support, $G$ is bounded and
\eqref{eq:regular-smooth-moments} already supplies such an integrable
$J$.

Under these smooth conditions, a consistent sandwich estimate is
\[
 \widehat H_n=P_nA(\cdot,\widehat u_n),\qquad
 \widehat\Sigma_n=n^{-1}\sum_i(s_i-\bar s_n)(s_i-\bar s_n)^{\mathrm T},
 \qquad \widehat V_n=\widehat H_n^{-1}\widehat\Sigma_n\widehat H_n^{-1},
\]
where $s_i=s(Y_i,\widehat u_n)$ and $\bar s_n=n^{-1}\sum_i s_i$.
Set $\widehat V_n=I_d$ when an inverse is undefined, an estimated
matrix is not positive definite, or $\widehat u_n$ uses its exceptional
definition; these events have probability tending to zero.
This Hessian estimator is used only under the smooth conditions.
For the circle, \cref{sec:supp-circle-model} instead estimates
curvature using the antipodal density.

\begin{proof}[\proofheading{\mainCref{prop:wald-recovery}}]
\noindent\emph{Localization.}
Let $P_n$ denote the empirical law. By Varadarajan's theorem and the
strong law, $P_n$ converges weakly to $P$ with convergence of second
moments about $o$ a.s.. Thus the quadratic Wasserstein distance
satisfies $W_2(P_n,P)\to0$ a.s.
\citep[Definition~6.8 and Theorem~6.9]{Villani2009a}.
Since $F_R(q)^{1/2}=W_2(R,\delta_q)$ for $R=P_n,P$, where $\delta_q$
is the point mass at $q$, the reverse triangle inequality and taking infima give
\begin{equation}
 \begin{aligned}
 \sup_{q\in\cQ}|F_n(q)^{1/2}-F_P(q)^{1/2}|
 &\le W_2(P_n,P)\to0,\\
 |(F_n^\star)^{1/2}-(F_P^\star)^{1/2}|
 &\le W_2(P_n,P)\to0
 \qquad\text{a.s.}.
 \end{aligned}
 \label{eq:regular-root-risk-law}
\end{equation}
For fixed $K\ge0$, define
\[
 L_n(K)\coloneqq
 \{q\in\cQ:F_n(q)\le F_n(m)+Kn^{-1/2}\wh\sigma_n(q,m)\}.
\]
By \eqref{eq:loss-risk-scale} and
$d(q,m)\le F_n(q)^{1/2}+F_n(m)^{1/2}$, every $q\in L_n(K)$ satisfies
\[
 F_n(q)^{1/2}\le(1+2Kn^{-1/2})F_n(m)^{1/2}.
\]
Indeed, for $x=F_n(q)^{1/2}$ and $y=F_n(m)^{1/2}$,
the defining inequality gives $(x-y)(x+y)\le2Kn^{-1/2}y(x+y)$;
cancel $x+y$ when positive, while $x+y=0$ is immediate.
Writing $\Delta_n=W_2(P_n,P)$, \eqref{eq:regular-root-risk-law} yields
\[
 0\leq \sup_{q\in L_n(K)}\{F_P(q)^{1/2}-(F_P^\star)^{1/2}\}
 \le2Kn^{-1/2}(F_P^\star)^{1/2}
       +(2+2Kn^{-1/2})\Delta_n\longrightarrow0\quad\text{a.s.}.
\]
By \eqref{eq:regular-risk-isolation}, $L_n(K)$ is therefore contained
in any fixed neighbourhood of $m$ with probability tending to one.

The selected comparisons are localized similarly. An anchor satisfies
$F_n(q)\le(1+b_n n^{-1/2})F_n^\star$, while $T_n(q,r)\ge-\tau_n$
and \eqref{eq:loss-risk-scale} imply
\[
 F_n(r)^{1/2}\le(1+2\tau_n n^{-1/2})F_n(q)^{1/2}.
\]
Since $b_n,\tau_n=o(n^{1/2})$ by \mainequationref{eq:adaptive-tuning-main},
\eqrefs{eq:regular-root-risk-law,eq:regular-risk-isolation} place both
endpoints of every selected pair in any fixed chart neighbourhood
with probability tending to one.

\medskip\noindent\emph{Local quadratic and variance approximations.}
Under (W3), apply the fixed-$P$ arguments of
\cref{lem:uniform-squares,lem:relative-variance} to
$\mathcal H_{\mathrm{loc}}$ to obtain
\[
 K_n^{\mathrm{loc}}\coloneqq\sup_{h\in\mathcal H_{\mathrm{loc}}}
 |\mathbb G_nh|=O_p(1),\qquad
 \eta_n\coloneqq\sup_{(q,r)\in\mathcal I_{\mathrm{loc}}}
 \left|\frac{\wh\sigma_n(q,r)}{\sigma_P(q,r)}-1\right|=o_p(1),
\]
where $\mathcal I_{\mathrm{loc}}$ indexes the positive-variance pairs
in $\phi(B(0,\delta_0))^2$.
Population zero-variance pairs obey \cref{lem:zero-variance}.
By (W1)--(W2) and positive definiteness of $\Sigma$ and $H$, there are
$c,C>0$ and $0<\delta<\delta_0$ such that, for $\|u\|\le\delta$,
\[
 F_P(\phi(u))-F_P^\star\ge c\|u\|^2,\qquad
 \sigma_P(\phi(u),m)\le C\|u\|,
\]
and $\sigma_P(\phi(u),m)>0$ whenever $u\ne0$.
On the localization event $L_n(K)\subseteq\phi(B(0,\delta))$,
every $q=\phi(u)\in L_n(K)$ with $u\ne0$ satisfies
\[
 \begin{aligned}
 c\|u\|^2&\le F_P(q)-F_P^\star\\
 &\le n^{-1/2}\{K\wh\sigma_n(q,m)
                +K_n^{\mathrm{loc}}\sigma_P(q,m)\}\\
 &\le Cn^{-1/2}\{K(1+\eta_n)+K_n^{\mathrm{loc}}\}\|u\|.
 \end{aligned}
\]
The second inequality uses the defining inequality for $L_n(K)$ and
$F_n(q)-F_n(m)=F_P(q)-F_P^\star+n^{-1/2}\sigma_P(q,m)\mathbb G_nh_{q,m}$.
The case $u=0$ is immediate. Thus
\begin{equation}
 \sup_{q\in L_n(K)}n^{1/2}\|\phi^{-1}(q)\|=O_p(1).
 \label{eq:regular-root-localization}
\end{equation}
The supremum is set to zero off the localization event.

Finally, $L_n(0)$ is contained with probability tending to one in the
image of a fixed compact coordinate ball. Since $F_n$ is continuous
and every point outside this image has empirical risk greater than
$F_n(m)$, a global empirical minimizer exists on this event.
Minimizing over the fixed compact set gives measurable choices that
are global minimizers on the same event.

The finite-dimensional class $\mathcal E$ is Donsker, so (W3)
and the joint finite-dimensional CLT make
$\mathcal H_{\mathrm{loc}}\cup\mathcal E$ Donsker.
The empirical and conditional multiplier processes are therefore
asymptotically equicontinuous on this union in $L^2(P)$.
Write $e_a=a^{\mathrm T}\psi/(a^{\mathrm T}\Sigma a)^{1/2}$ for
$\|a\|=1$. For $a=u/\|u\|$ and $0<\|u\|\le\delta$, (W1) gives
\[
 \begin{aligned}
 \left|\frac{\sigma_P(m,\phi(u))}{\|u\|}
               -(a^{\mathrm T}\Sigma a)^{1/2}\right|
 &\le\epsilon_A(\delta),\\
 \|h_{m,\phi(u)}-e_a\|_2
 &\le\frac{2\epsilon_A(\delta)}{\lambda_{\min}(\Sigma)^{1/2}}.
 \end{aligned}
\]
Decompose the normalized empirical remainder into
$\{\sigma_P(m,\phi(u))/\|u\|\}\mathbb G_n(h_{m,\phi(u)}-e_a)$
and the difference of the two scales times $\mathbb G_ne_a$.
For fixed $R<\infty$, equicontinuity and
$\sup_{e\in\mathcal E}|\mathbb G_ne|=O_p(1)$ give
\[
 \sup_{0<\|u\|\le Rn^{-1/2}}
 \frac{|\mathbb G_n(\ell_{m,\phi(u)}-u^{\mathrm T}\psi)|}{\|u\|}
 =o_p(1).
\]
Together with (W2), this yields
\begin{equation}
 \begin{aligned}
 &\sup_{\|z\|\le R}
 \left|n\{f_n(z/n^{1/2})-f_n(0)\}
              -\tfrac12z^{\mathrm T}Hz+z^{\mathrm T}W_n\right|=o_p(1),\\
 &\sup_{\|z\|,\|w\|\le R}
 \left|n^{1/2}\wh\sigma_n(\phi(z/n^{1/2}),\phi(w/n^{1/2}))
              -\|\Sigma^{1/2}(z-w)\|\right|=o_p(1).
 \end{aligned}
 \label{eq:regular-quadratic-approximation}
\end{equation}
For the second assertion in \eqref{eq:regular-quadratic-approximation}, the identity
$\ell_{\phi(u),\phi(v)}=\ell_{m,\phi(v)}-\ell_{m,\phi(u)}$
and (W1) imply
\[
 \left|\sigma_P(\phi(u),\phi(v))
              -\|\Sigma^{1/2}(v-u)\|\right|
 \le\epsilon_A(\delta)(\|u\|+\|v\|).
\]
Apply uniform local studentization; population zero-variance pairs
satisfy the same absolute-error bound by \cref{lem:zero-variance}.

All empirical minimizers have root-$n$ tight coordinates by
\eqref{eq:regular-root-localization}. On a fixed ball containing
these coordinates and $H^{-1}W_n$, their minimizing inequalities give
\[
 \tfrac12\lambda_{\min}(H)
       \|n^{1/2}\widehat u_n-H^{-1}W_n\|^2\le2\sup_{\|z\|\le R}
 \left|n\{f_n(z/n^{1/2})-f_n(0)\}
              -\tfrac12z^{\mathrm T}Hz+z^{\mathrm T}W_n\right|=o_p(1).
 \]

Such a ball contains all empirical minimizers and $H^{-1}W_n$
with arbitrarily high probability for sufficiently large $R$.
Consequently, simultaneously for all empirical minimizers,
\begin{equation}
 n^{1/2}\widehat u_n=H^{-1}W_n+o_p(1),\qquad
 W_n\rightsquigarrow N_d(0,\Sigma).
 \label{eq:regular-mean-expansion}
\end{equation}

\medskip\noindent\emph{Critical value and coverage.}
Localization of selected endpoints, uniform local studentization,
\eqref{eq:regular-direction-stability} and conditional multiplier
equicontinuity give
\[
 S_n^{\#,\mathrm L}\le\sup_{e\in\mathcal E}\mathbb G_n^\#e
       +o_{p_\xi}(1)
 =\|\Sigma^{-1/2}\mathbb G_n^\#\psi\|+o_{p_\xi}(1)
\]
in probability. In particular, $c_n^{\mathrm L}=O_p(1)$.
For the reverse bound,
\eqrefs{eq:regular-quadratic-approximation,eq:regular-mean-expansion} give
\[
 n\{F_n(m)-F_n^\star\}
 =\tfrac12W_n^{\mathrm T}H^{-1}W_n+o_p(1)=O_p(1).
\]
Moreover, $F_P^\star>0$, since otherwise $P$ is a point mass at $m$
and (W1) contradicts $\Sigma>0$.
Equation \eqref{eq:regular-root-risk-law} gives
$F_n^\star\to_p F_P^\star$, so $b_n\to\infty$ puts
$m\in\wh M_n$ with probability tending to one.
Uniformly for $\|a\|=1$, (W1) and (W2) give
\[
 \|h_{m,\phi(a/n)}-e_a\|_2=o(1),\qquad
 n^{1/2}\lambda_P(m,\phi(a/n))=O(n^{-1/2}).
\]
Local studentization and empirical equicontinuity imply
\[
 \sup_{\|a\|=1}|T_n(m,\phi(a/n))-\mathbb G_ne_a|=o_p(1).
\]
These comparisons have positive empirical variance with probability
tending to one. Since $\sup_{e\in\mathcal E}|\mathbb G_ne|=O_p(1)$
and $\tau_n\to\infty$, the selector retains all of them on the same
event. Conditional multiplier equicontinuity gives the matching
lower bound. Thus, conditionally in probability,
\begin{equation}
 S_n^{\#,\mathrm L}
 =\|\Sigma^{-1/2}\mathbb G_n^\#\psi\|+o_{p_\xi}(1),\qquad
 c_n^{\mathrm L}=c_\alpha+o_p(1).
 \label{eq:regular-critical-value}
\end{equation}
The conditional covariance of
$\Sigma^{-1/2}\mathbb G_n^\#\psi$ converges to $I_d$ by the law of
large numbers. Continuity and strict increase of the $\chi_d$
distribution at $c_\alpha$ justify the quantile conclusion.

For coverage, $T_n(m,r)>0$ implies $F_n(r)<F_n(m)$, so
$r\in L_n(0)$. Thus all positive statistics have local comparators
by \eqref{eq:regular-root-localization}.
Dropping their nonpositive population drift and using (W1),
local studentization and empirical equicontinuity bounds their
supremum by $\sup_{e\in\mathcal E}\mathbb G_ne+o_p(1)$.
Nonpositive statistics obey the same upper bound, since this
supremum is nonnegative. The comparisons $(m,\phi(a/n))$ give the
reverse bound. Hence
\begin{equation}
 \sup_{r:\wh\sigma_n(m,r)>0}T_n(m,r)
 =\|\Sigma^{-1/2}W_n\|+o_p(1).
 \label{eq:regular-mean-statistic}
\end{equation}
With probability tending to one, no zero-empirical-variance comparison
rejects $m$: rejection would put $r\ne m$ in $L_n(0)$,
where (W1) and local studentization ensure positive empirical
variance. Population zero-variance comparisons also satisfy inversion
by \cref{lem:zero-variance}.
More generally, if a nonnegative common cutoff satisfies
$c_n\to_p c$ for a constant $c\ge0$, then
\[
 \pr_P(m\in\mathcal C_n(c_n))\longrightarrow\pr(\chi_d\le c).
\]
Indeed, apart from an event of probability tending to zero, membership
is equivalent to the statistic in \eqref{eq:regular-mean-statistic}
being at most $c_n$; the CLT and continuity of the
$\chi_d$ distribution give the limit. Hence a constant limiting cutoff
smaller than $c_\alpha$ cannot attain nominal asymptotic coverage within
the inversion family.
Equations
\eqrefs{eq:regular-critical-value,eq:regular-mean-statistic}
and the CLT give $\pr_P(m\in C_n^{\mathrm L})\to1-\alpha$.

\medskip\noindent\emph{Wald sets from penalized minimization.}
The empirical standard deviation $\wh\sigma_n$ is a pseudometric,
being the empirical $L^2$ distance between centred losses.
For a target $q_0$ and $c\ge0$, let
\[
 \widetilde q\in\argmin_q
 \{F_n(q)+cn^{-1/2}\wh\sigma_n(q,q_0)\}.
\]
The minimizing inequality and the triangle inequality give, for every $r$,
\begin{equation}
 F_n(\widetilde q)-F_n(r)
 \le cn^{-1/2}\{\wh\sigma_n(r,q_0)-\wh\sigma_n(\widetilde q,q_0)\}
 \le cn^{-1/2}\wh\sigma_n(\widetilde q,r).
 \label{eq:regular-penalized-inclusion}
\end{equation}
Thus $\widetilde q\in\mathcal C_n(c)$.
For fixed $K$, these minima exist simultaneously for $0\le c\le K$
and all targets with probability tending to one: the sublevel below
the objective's value at $m$ lies in $L_n(K)$ by the reverse triangle
inequality, so continuity reduces minimization to a fixed compact
chart ball. Comparing with $m$ also gives
$\mathcal C_n(c)\subseteq L_n(K)$ for $0\le c\le K$.
Thus all accepted points and all penalized minimizers satisfy
\eqref{eq:regular-root-localization}, uniformly over targets and cutoffs.

In root-$n$ coordinates put
\[
 A_n(c)=\{z:\|\Sigma^{-1/2}(Hz-W_n)\|\le c\}.
\]
For a target $z$, the function
$w\mapsto\tfrac12w^{\mathrm T}Hw-w^{\mathrm T}W_n
                  +c\|\Sigma^{1/2}(w-z)\|$
is strongly convex with curvature $\lambda_{\min}(H)$.
Its unique minimizer $v$ satisfies
$Hv-W_n=-c\Sigma^{1/2}a$ for some $\|a\|\le1$, so $v\in A_n(c)$;
moreover, $v=z$ if and only if $z\in A_n(c)$.
By \eqref{eq:regular-quadratic-approximation}, for fixed $R,K$,
\[
 \begin{aligned}
 &n\{f_n(w/n^{1/2})-f_n(0)\}
   +cn^{1/2}\wh\sigma_n(\phi(w/n^{1/2}),\phi(z/n^{1/2}))\\
 &\quad=\tfrac12w^{\mathrm T}Hw-w^{\mathrm T}W_n
       +c\|\Sigma^{1/2}(w-z)\|+o_p(1),
 \end{aligned}
\]
where the remainder is uniform over $\|z\|,\|w\|\le R$ and $0\le c\le K$.

If $z\in n^{1/2}\phi^{-1}(\mathcal C_n(c))$, acceptance says that
the empirical penalized objective at $z$ is no larger than at $v$.
When $z,v$ lie in the ball, strong convexity and the uniform
approximation therefore give
\begin{equation}
 \begin{aligned}
 \tfrac12\lambda_{\min}(H)\|z-v\|^2
 &\le\tfrac12z^{\mathrm T}Hz-z^{\mathrm T}W_n
       -\tfrac12v^{\mathrm T}Hv+v^{\mathrm T}W_n
       -c\|\Sigma^{1/2}(v-z)\|\le o_p(1),\\
 \operatorname{dist}(z,A_n(c))&\le\|z-v\|=o_p(1).
 \end{aligned}
 \label{eq:regular-wald-outer-bound}
\end{equation}
Conversely, take $z\in A_n(c)$ and an empirical penalized minimizer
with target $\phi(z/n^{1/2})$, whose root-$n$ coordinate is
$\widetilde z$. Its image is accepted by
\eqref{eq:regular-penalized-inclusion}, and $z$ minimizes the
quadratic penalized objective. When $z,\widetilde z$ lie in the ball,
the same argument gives
\begin{equation}
 \begin{aligned}
 \tfrac12\lambda_{\min}(H)\|\widetilde z-z\|^2
 &\le\tfrac12\widetilde z^{\mathrm T}H\widetilde z
       -\widetilde z^{\mathrm T}W_n
       -\tfrac12z^{\mathrm T}Hz+z^{\mathrm T}W_n
       +c\|\Sigma^{1/2}(\widetilde z-z)\|\le o_p(1),\\
 \operatorname{dist}\bigl(z,n^{1/2}\phi^{-1}(\mathcal C_n(c))\bigr)
 &\le\|\widetilde z-z\|=o_p(1).
 \end{aligned}
 \label{eq:regular-wald-inner-bound}
\end{equation}
The bounds in \eqrefs{eq:regular-wald-outer-bound,eq:regular-wald-inner-bound}
use a common $o_p(1)$ upper bound, uniform over all the
indicated points and $0\le c\le K$ on each fixed ball.
All points of $A_n(c)$ and all quadratic penalized minimizers are
uniformly tight for $c\le K$, since
$A_n(c)=H^{-1}W_n+H^{-1}\Sigma^{1/2}
\{a\in\mathbb R^d:\|a\|\le c\}$.
Accepted points and empirical penalized minimizers obey the same
control by \eqref{eq:regular-root-localization}.
Letting $R\to\infty$ and then $K\to\infty$ yields, for every
nonnegative tight cutoff $c_n$,
\begin{equation}
 d_H\bigl(n^{1/2}\phi^{-1}(\mathcal C_n(c_n)),A_n(c_n)\bigr)=o_p(1).
 \label{eq:regular-set-approximation}
\end{equation}
By \eqrefs{eq:regular-mean-expansion,eq:regular-critical-value}
and the assumed consistency of $\widehat V_n$, the Wald ellipsoid
lies inside the chart with probability tending to one, and
\[
 \begin{aligned}
 &d_H\bigl(A_n(c_n^{\mathrm L}),
                 n^{1/2}\phi^{-1}(\mathcal E_n(c_\alpha))\bigr)\\
 &\quad\le\|H^{-1}W_n-n^{1/2}\widehat u_n\|
       +c_n^{\mathrm L}\|V^{1/2}-\widehat V_n^{1/2}\|_{\mathrm{op}}+|c_n^{\mathrm L}-c_\alpha|
                   \|\widehat V_n^{1/2}\|_{\mathrm{op}}=o_p(1).
 \end{aligned}
\]
The forward Lipschitz chart proves the asserted metric Hausdorff
conclusion. If \maincref{ass:main-process} also holds for
$\mathcal P=\{P\}$, the global calibration theorem and the same argument apply with
$c_n^G\to_p c_{P,\alpha}^G$;
$\mathcal E\subseteq\overline{\mathcal H_P}^{\,L^2(P)}$ by (W1)
implies $c_{P,\alpha}^G\ge c_\alpha$.

\medskip\noindent\emph{Smooth alternative.}
The square-integrable score bound and differentiation at the
minimum give $Ps(Y,0)=0$.
Under \eqref{eq:regular-smooth-moments}, score continuity and
dominated convergence give
\[
 P\sup_{\|u\|\le\delta}\|s(Y,u)-s(Y,0)\|^2\to0,
\]
and the Hessian bound gives
\[
 \sup_{\|u\|\le\delta}\|P_nA(\cdot,u)-H\|_{\mathrm{op}}
 \le\|P_nA(\cdot,0)-H\|_{\mathrm{op}}+P_nJ_\delta.
\]
For $\|u\|,\|v\|\le\delta$, the integral increment formula implies
\[
 |g(y,u)-g(y,v)-(u-v)^{\mathrm T}s(y,0)|
 \le\|u-v\|\sup_{\|w\|\le\delta}\|s(y,w)-s(y,0)\|.
\]
Its centred empirical $L^2$ norm is at most
$\|u-v\|\{P_n\sup_{\|w\|\le\delta}
\|s(Y,w)-s(Y,0)\|^2\}^{1/2}$.
The law of large numbers, first as $n\to\infty$ and then as
$\delta\downarrow0$, gives relative local variance approximation
and positive empirical curvature. On the localization event,
Taylor's formula bounds $f_n(u)-f_n(0)$ below by
$c\|u\|^2-\|P_ns(\cdot,0)\|\|u\|$ for a fixed $c>0$.
Together with $\wh\sigma_n(\phi(u),m)\le C\|u\|$ and the score CLT,
this proves \eqref{eq:regular-root-localization}.
Taylor's formula then gives
\eqrefs{eq:regular-quadratic-approximation,eq:regular-mean-expansion}
with $\psi=-s(\cdot,0)$, without a Donsker assumption.
The preceding penalized-minimization argument proves the set
approximation \eqref{eq:regular-set-approximation}.

The integral score formula and \eqref{eq:regular-score-multiplier}
approximate all selected normalized multiplier contrasts by score
directions. The comparisons $(m,\phi(a/n))$ again give the matching
lower bound in \eqref{eq:regular-critical-value}.
For coverage, positive empirical curvature bounds the statistic at
$m$ above by $\|\Sigma^{-1/2}\mathbb G_n\psi\|+o_p(1)$, and
arbitrarily short coordinate comparisons attain this bound.
Any comparator with $F_n(r)<F_n(m)$ belongs to $L_n(0)$, so the
localization argument excludes remote comparators, including
zero-empirical-variance rejections.
This proves \eqref{eq:regular-mean-statistic} and coverage.
Finally, the same $L^1$ Hessian bound and $L^2$ score continuity give
$\widehat H_n\to_p H$ and $\widehat\Sigma_n\to_p\Sigma$, validating
the stated sandwich estimate.
\end{proof}

For the circular family in \mainequationref{eq:circle-family}, fixed
$\lambda>0$ gives a unique isolated mean at zero by
\eqref{eq:circle-family-risk} and compactness. With the angular chart
at zero, \eqref{eq:circle-strict-secant} gives (W1) and the direction
condition in (W3), while \eqref{eq:circle-family-risk} gives (W2) with
$H=\lambda$. The empirical-process argument in the proof of
\maincref{cor:circle-exactness-main} supplies the Donsker part of (W3).
Thus the regular local conditions of \maincref{prop:wald-recovery} hold.

\subsection{Confidence radii}
\label{sec:supp-radii}
\begin{proof}[\proofheading{\mainCref{prop:adaptive-radius}}]

Let $E_{n,P}$ be the common support event of
\cref{lem:zero-variance} intersected with the defining event of
\cref{lem:uniform-quantiles}, so
$\sup_P\pr_P(E_{n,P}^c)\to0$. Work on $E_{n,P}$.

Inversion and the empirical-process bound give, simultaneously for
$q\in\mathcal C_n(c_n)$ and $m\in M_P$,
\begin{equation}
 \begin{aligned}
 n^{1/2}\{F_P(q)-F_P^\star\}
 &\le c_n\wh\sigma_n(q,m)+\sigma_P(q,m)\|\mathbb G_n\|_{\mathcal H_P}\\
 &\le\sigma_P(q,m)
 \left\{c_n\sup_{(r,s)\in\mathcal I_P}
            \frac{\wh\sigma_n(r,s)}{\sigma_P(r,s)}+\|\mathbb G_n\|_{\mathcal H_P}\right\}.
 \end{aligned}
 \label{eq:proof-radius-master-bound}
\end{equation}
Here \cref{lem:zero-variance} covers zero-variance comparisons.
Let $A_n$ be the quantity in braces in
\eqref{eq:proof-radius-master-bound}, divided by $n^{1/2}$, and set
$A_n=0$ off this event. By \cref{lem:uniform-squares,lem:relative-variance}
and $c_n=O_{\mathcal P}(1)$,
$A_n=O_{\mathcal P}(n^{-1/2})$.
Taking the infimum over $m\in M_P$ and using
\eqref{eq:loss-risk-scale} give
\[
 F_P(q)-F_P^\star
 \le A_n\inf_{m\in M_P}\sigma_P(q,m)
 \le2A_n\{F_P(q)+F_P^\star\}.
\]
Thus, on $E_{n,P}\cap\{A_n\le1/4\}$,
\[
 \sup_{q\in\mathcal C_n(c_n)}\{F_P(q)-F_P^\star\}
 \le8F_P^\star A_n\le8\overline F A_n.
\]

For $\epsilon>0$, put
\[
 \kappa(\epsilon)\coloneqq
 \inf_{P\in\mathcal P}\inf_{q:\ d(q,M_P)\ge\epsilon}
 \{F_P(q)-F_P^\star\}>0.
\]
On the preceding event, $\operatorname{rad}(\mathcal C_n(c_n);M_P)>\epsilon$
implies $8\overline F A_n\ge\kappa(\epsilon)$. Hence
\[
 \begin{aligned}
 &\sup_{P\in\mathcal P}\pr_P
 (\operatorname{rad}(\mathcal C_n(c_n);M_P)>\epsilon)\\
 &\quad\le\sup_P\pr_P(E_{n,P}^c)+\sup_P\pr_P(A_n>1/4)+\sup_P\pr_P(8\overline F A_n\ge\kappa(\epsilon))\to0.
 \end{aligned}
\]
Thus $\operatorname{rad}(\mathcal C_n(c_n);M_P)=o_{\mathcal P}(1)$.

Consequently, with probability tending to one
uniformly over $P$, every accepted non-mean
satisfies $0<d(q,M_P)\le r_0$.
The margin condition \mainequationref{eq:margin-scale-main} and
\cref{lem:automatic-linear-scale} then give
\[
 a\,d(q,M_P)^\beta
 \le F_P(q)-F_P^\star
 \le A_n\inf_{m\in M_P}\sigma_P(q,m)
 \le2\overline F^{1/2}A_n\,d(q,M_P).
\]
Dividing by $d(q,M_P)>0$ and taking the supremum over accepted points gives
\[
 \operatorname{rad}(\mathcal C_n(c_n);M_P)
 \le\left\{2a^{-1} \overline F^{1/2}A_n\right\}^{1/(\beta-1)}
 =O_{\mathcal P}\bigl(n^{-1/\{2(\beta-1)\}}\bigr),
\]
proving \mainequationref{eq:adaptive-radius-main}.
Finally, $c_n^G=O_{\mathcal P}(1)$ by
\eqref{eq:uniform-pregaussian} and \maincref{thm:calibration}, and
$0\le c_n^{\mathrm L}\le c_n^G$. Thus the result applies to both calibrations.
\end{proof}
\section{Pointwise refinements and limits of uniform validity}
\label{sec:pointwise-refinements}
\subsection{Scale-restricted calibration}
\label{sec:supp-scale-refinement}
For $v\ge0$, restrict the calibration selector to
\[
 \begin{aligned}
 &\widehat{\mathcal B}_n^{\mathrm{ref}}(v)
 \coloneqq\{(q,r)\in\widehat{\mathcal B}_n^{\mathrm L}:
                         \wh\sigma_n(q,r)\ge v\},\qquad
 S_n^{\#,\mathrm{ref}}(v)\coloneqq
 \left(\sup_{(q,r)\in\widehat{\mathcal B}_n^{\mathrm{ref}}(v)}T_n^\#(q,r)\right)_+,\\
 &c_n^{\mathrm{ref}}(v)\coloneqq
 q_{1-\alpha}(S_n^{\#,\mathrm{ref}}(v)\mid Y_1,\ldots,Y_n),\qquad
 C_n^{\mathrm{ref}}(v)\coloneqq\mathcal C_n(c_n^{\mathrm{ref}}(v)).
 \end{aligned}
\]
Samplewise,
$c_n^{\mathrm{ref}}(v)\le c_n^{\mathrm L}$ and
$C_n^{\mathrm{ref}}(v)\subseteq C_n^{\mathrm L}\subseteq C_n^G$;
at $v=0$ the refined and localized procedures coincide.
The corresponding population directions are
$\mathcal H_{P,n}^{\mathrm{ref}}(v)\coloneqq
\{h_{q,r}:(q,r)\in\widehat{\mathcal B}_n^{\mathrm{ref}}(v)\cap\mathcal I_P\}$.
Define also
\[
 \mathcal A_P(\eta,v)\coloneqq
 \{h_{m,r}:m\in M_P,\ (m,r)\in\mathcal I_P,\
 \lambda_P(m,r)\le\eta,\ \sigma_P(m,r)\ge v\}.
\]
The anchor modulus is
$u_P(\delta)\coloneqq\sup_{q\in M_P(\delta)}
\inf_{m\in M_P}\sigma_P(q,m)$. Put
$\Delta_n\coloneqq2b_nF_n^\star/n^{1/2}$.

For two defined normalized contrasts,
\begin{equation}
 \|h_{q,r}-h_{m,r}\|_2
 \leq 2\frac{\sigma_P(q,m)}{\sigma_P(q,r)},\qquad
 |\sigma_P(q,r)-\sigma_P(m,r)|\leq\sigma_P(q,m).
 \label{eq:anchor-replacement}
\end{equation}
These follow from $\ell_{q,r}-\ell_{m,r}=\ell_{q,m}$, the reverse
triangle inequality and normalization.
\begin{theorem}[Pointwise exact calibration with a scale cutoff]
\label{thm:scale}
Fix $P$ and $0<\alpha<1/2$. Assume \maincref{ass:main-process}
for $\mathcal P=\{P\}$, \mainequationref{eq:adaptive-tuning-main},
and $\mathcal H_{P,0}\ne\varnothing$.
Let $v_n\ge0$ be measurable with $v_n\to_p0$. Suppose
\begin{equation}
 \frac{u_P(\Delta_n)
               +\Delta_n}{v_n}\to_p0
 \label{eq:scale-separation-main}
\end{equation}
and
\begin{equation}
 \sup_{h\in\mathcal H_{P,0}}
 \inf_{g\in\mathcal A_P(\tau_n/(16n^{1/2}),2v_n)}
 \rho_P(h,g)\to_p0.
 \label{eq:scale-retention-main}
\end{equation}
The ratio in \eqref{eq:scale-separation-main} is $+\infty$ on
$\{v_n=0\}$. Then
$d_{H,\rho_P}(\overline{\mathcal H_{P,n}^{\mathrm{ref}}(v_n)}^{\,L^2(P)},
\mathcal H_{P,0})\to_p0$,
$c_n^{\mathrm{ref}}(v_n)\to_pc_{P,\alpha}^0$ and
$\pr_P(M_P\subseteq C_n^{\mathrm{ref}}(v_n))\to1-\alpha$.
\end{theorem}
\begin{proof}
Work on \eqrefs{eq:proof-gms-event,eq:relaxed-tube-inclusions}
and $\{v_n>0\}$. For a selected pair $(q,r)$,
$q\in M_P(\Delta_n)$ and $\sigma_P(q,r)\ge v_n/2$.
Choose $m\in M_P$ with
$\sigma_P(q,m)\le u_P(\Delta_n)+v_n/n$.
By \eqref{eq:scale-separation-main}, the reverse triangle inequality and
\eqref{eq:anchor-replacement},
with probability tending to one, simultaneously for all selected pairs,
\begin{equation}
 \sigma_P(m,r)\ge\tfrac12\sigma_P(q,r)\ge v_n/4>0,\qquad
 \|h_{q,r}-h_{m,r}\|_2
 \le4\{u_P(\Delta_n)/v_n+n^{-1}\}.
 \label{eq:scale-anchor-error}
\end{equation}
The slack-selection inequality gives
$F_P(r)-F_P(q)\le3\tau_n\sigma_P(q,r)/n^{1/2}$.
Together with $F_P(q)-F_P^\star\le\Delta_n$, this yields
\begin{equation}
 0\le\lambda_P(m,r)
 \le\frac{6\tau_n}{n^{1/2}}
       +\frac{4\Delta_n}{v_n}
 \to_p0.
 \label{eq:scale-slack-error}
\end{equation}
Both error bounds vanish in probability. Thus \cref{lem:compact-binding}
and \eqrefs{eq:scale-anchor-error,eq:scale-slack-error} imply that the
Hausdorff excess of the selected directions and their closure over
$\mathcal H_{P,0}$ is $o_p(1)$.

Conversely, each pair generating
$\mathcal A_P(\tau_n/(16n^{1/2}),2v_n)$ has its
first endpoint in $\widehat M_n$, empirical standard deviation at
least $v_n$, and $T_n(m,r)\ge-3\tau_n/16>-\tau_n$ on
\eqref{eq:proof-gms-event}. Thus
\[
 \mathcal A_P(\tau_n/(16n^{1/2}),2v_n)
 \subseteq\mathcal H_{P,n}^{\mathrm{ref}}(v_n).
\]
By \eqref{eq:scale-retention-main}, the reverse Hausdorff excess is
$o_p(1)$ and the selected class is nonempty with probability tending
to one.
\Cref{lem:uniform-squares,lem:relative-variance} make the multiplier
studentization error conditionally negligible.
\Cref{lem:moving-sets,lem:gaussian-quantile} then give
$c_n^{\mathrm{ref}}(v_n)\to_pc_{P,\alpha}^0$;
\maincref{thm:binding-main} and \eqref{eq:coverage-event-on-good-set}
give exact coverage.
\end{proof}

\subsection{A sticky transition}
Let $\cQ$ be the tripod formed from three copies of $[0,\infty)$,
joined at $o$, and write $(j,t)$ for distance $t$ on branch $j$.
Observations have distance one and branch probabilities
\[
 (p_1,p_2,p_3)=((1+\theta)/2,(1-\theta)/4,(1-\theta)/4),
 \qquad 0\le\theta\le\theta_*,\quad 0<\theta_*<1/2.
\]
The three support probabilities have a common positive lower bound,
so \cref{lem:ep-sufficient} verifies \maincref{ass:main-process}
uniformly on this family. Let $\pr_\theta,\E_\theta$ refer to $P_\theta$,
and put $X_j=2\mathds{1}\{Y=(j,1)\}-1$ and $X=X_1$.
Directly,
\begin{equation}
 F_{P_\theta}((j,t))=1+2(1-2p_j)t+t^2,\quad
 \ell_{o,(j,t)}=2t\{X_j-\E_\theta X_j\},\quad
 \sigma_{P_\theta}(o,(j,t))=4t\{p_j(1-p_j)\}^{1/2}.
 \label{eq:general-tripod-formulas}
\end{equation}
In particular, $\E_\theta X=\theta$, $\Var_\theta(X)=1-\theta^2$,
$F_{P_\theta}((1,t))=1-2\theta t+t^2$, and
$M_{P_0}=\{o\}$, $M_{P_\theta}=\{(1,\theta)\}$ for $\theta>0$.
At $P_0$, $F_{P_0}^\star=1$ and, for $t>0$, the three branch slacks are
$t/2$, $(1+t)/3^{1/2}$ and $(1+t)/3^{1/2}$, with $h_{o,(1,t)}=X$.
Consequently,
\[
 \mathcal H_{P_0,0}=\{X\},\qquad T_{P_0}^0\sim N(0,1),\qquad
 c_{P_0,\alpha}^0=z_{1-\alpha},
\]
where $z_p=\Phi^{-1}(p)$.

For this refinement example, take, for $n\ge3$,
\[
 b_n=(\log\log n)^{1/2},\qquad \tau_n=n^{3/8},\qquad
 v_n=F_n^\star n^{-1/5}.
\]
This $\tau_n$ differs from the numerical Local tuning $n^{1/5}$.
Since $F_n^\star\to_{P_0}1$, $v_n\to_{P_0}0$ and $v_n>0$ with
probability tending to one. On $\{v_n>0\}$,
\[
 \frac{\Delta_n^{1/2}}{v_n}
 =\frac{(2b_n)^{1/2}}{(F_n^\star)^{1/2}}n^{-1/20}\to_{P_0}0,\qquad
 \frac{n^{1/2}v_n}{\tau_n}
 =F_n^\star n^{-3/40}\to_{P_0}0.
 \]
The risk and scale formulas in \eqref{eq:general-tripod-formulas} give
$u_{P_0}(\delta)\le2\delta^{1/2}$, verifying
\eqref{eq:scale-separation-main}.
The comparison $(o,(1,v_n))$ has direction $X$, scale $2v_n$ and
slack $v_n/2\le\tau_n/(16n^{1/2})$ with probability tending to one.
Thus \eqref{eq:scale-retention-main} holds, and \cref{thm:scale} yields
\[
 c_n^{\mathrm{ref}}(v_n)\to_{P_0}z_{1-\alpha},\qquad
 \pr_0(M_{P_0}\subseteq C_n^{\mathrm{ref}}(v_n))\to1-\alpha.
\]

Local retains both $(o,(1,n^{-1}))$ and its reverse with probability
tending to one: \cref{lem:bounded-ball-localization} retains both
endpoints, and their empirical statistics are $O_p(1)$.
Their conditional multiplier statistics are opposite standard normals, so
\[
 \pr_0(c_n^{\mathrm L}\ge z_{1-\alpha/2})\to1.
\]
The following bound applies to every common scalar cutoff in
\mainequationref{eq:inversion}.
\begin{proposition}[Nonuniformity of oracle refinement]
\label{prop:adaptation-barrier}
Let $c_n$ be any nonnegative measurable common scalar cutoff with
$c_n\to_{P_0}c_*\in[0,\infty)$. For every $a>0$ and
$\theta_n=a/n^{1/2}$,
\[
 \limsup_n\pr_{P_{\theta_n}}
 \left(M_{P_{\theta_n}}\subseteq\mathcal C_n(c_n)\right)
 \le2\Phi(c_*)-1.
\]
Consequently, uniform nominal coverage on this family requires
$c_*\ge z_{1-\alpha/2}$ whenever the boundary cutoff has a constant
limit. In particular, $c_*=z_{1-\alpha}$ yields the upper bound
$1-2\alpha<1-\alpha$.
\end{proposition}
\begin{proof}
For $\theta_n=a/n^{1/2}$, the product likelihood ratio relative to $P_0$ is
$L_n=\prod_i(1+\theta_nX_i)$ and
$\E_0[L_n^2]=(1+\theta_n^2)^n\le e^{a^2}$.
Cauchy--Schwarz gives
$\pr_{\theta_n}(A_n)\le e^{a^2/2}\{\pr_0(A_n)\}^{1/2}$
for every measurable event $A_n$. Applied to
$\{|c_n-c_*|>\epsilon\}$, this yields $c_n\to_{P_{\theta_n}}c_*$.
Set $\bar X_n=n^{-1}\sum_iX_i$ and
$\wh s_n^2=1-\bar X_n^2$. On $\{\wh s_n>0\}$, for
$m_{\theta_n}=(1,\theta_n)$ and $r=(1,\theta_n+u)$,
\[
 T_n(m_{\theta_n},r)=\operatorname{sign}(u)Z_n
       -\frac{n^{1/2}|u|}{2\wh s_n},\qquad
 Z_n=\frac{n^{1/2}(\bar X_n-\theta_n)}{\wh s_n}.
\]
Set $Z_n=0$ on $\{\wh s_n=0\}$.
For every sufficiently large $n$, $\theta_n>0$, so $u\to0$
from either sign is admissible. Thus whole-set coverage implies
$|Z_n|\le c_n$. The bounded triangular Bernoulli array satisfies
Lindeberg's condition, $1-\theta_n^2\to1$ and $\wh s_n\to_{P_{\theta_n}}1$;
hence $Z_n\rightsquigarrow N(0,1)$ and $\pr_{\theta_n}(\wh s_n=0)\to0$.
Continuity of $|Z|$ and cutoff convergence give
\[
 \limsup_n\pr_{\theta_n}(M_{P_{\theta_n}}\subseteq\mathcal C_n(c_n))
 \le\lim_n\pr_{\theta_n}(|Z_n|\le c_n)=2\Phi(c_*)-1.
\]
\end{proof}
Applying \cref{prop:adaptation-barrier} to the refinement above gives,
for every $a>0$ and $\theta_n=a/n^{1/2}$,
\[
 \limsup_n\pr_{\theta_n}
 (M_{P_{\theta_n}}\subseteq C_n^{\mathrm{ref}}(v_n))
 \le1-2\alpha<1-\alpha.
\]
Thus its exact coverage at $P_0$ does not extend to uniform nominal
coverage on this family.

\section{Simulation details}
\label{sec:supp-numerics}
\subsection{Ranking metrics and means}
\label{sec:supp-population}
Let $r_\sigma(i)$ denote the rank of object $i$ in ordering $\sigma$.
Euclidean and footrule distances are
\[
 d_E(\sigma,\pi)=\|r_\sigma-r_\pi\|_2,\qquad
 d_F(\sigma,\pi)=\sum_{i=1}^4|r_\sigma(i)-r_\pi(i)|.
\]
The Kendall distance is
\[
 d_K(\sigma,\pi)
 =\sum_{1\le i<j\le4}
 \mathds{1}\!\left\{
 [r_\sigma(i)-r_\sigma(j)]
 [r_\pi(i)-r_\pi(j)]<0
 \right\}.
\]
It equals the minimum number of adjacent swaps required to
transform one ordering into the other \citep{DiaconisGraham1977}. Rational enumeration gives $F_P^E(q)=15-2r_q(4)$ and
\[
 \begin{array}{c|c|c|c}
 \text{Metric}&M_P&F_P^\star&\text{Minimum nonmean risk gap}\\ \hline
 d_E&\{\sigma:r_\sigma(4)=4\}&7&2\\
 d_K&\{1234,2314,3124\}&83/12&1/2\\
 d_F&\{1324,2134,3214\}&56/3&4/3
 \end{array}
\]
Full support ensures positive contrast variance for distinct candidates;
the finite-support criterion in \cref{sec:supp-entropy} and positive
nonmean risk gaps give pointwise exact localized calibration.

Table~\ref{tab:ranking-cardinalities} reports finite-ranking confidence-set
cardinalities for the same paired datasets.

\begin{table}[htbp]
\centering
\caption[Finite-ranking confidence-set cardinalities]{Finite-ranking confidence-set cardinalities over 1000 paired datasets
per condition. Global and Local are mean cardinalities; Difference is
the mean $\#(C_n^G\setminus C_n^{\mathrm L})$. Reduced is the percentage
with $C_n^{\mathrm L}\subsetneq C_n^G$.}
\label{tab:ranking-cardinalities}
\begin{tabular}{@{}lrrrrr@{}}
\hline
Metric & $n$ & Global & Local & Difference & Reduced (\%) \\
\hline
Euclidean & 50 & 6.853 & 6.487 & 0.366 & 15.7 \\
 & 100 & 5.962 & 5.850 & 0.112 & 4.5 \\
 & 200 & 5.914 & 5.806 & 0.108 & 4.5 \\
\hline
Kendall & 50 & 10.624 & 8.730 & 1.894 & 81.8 \\
 & 100 & 7.667 & 6.460 & 1.207 & 66.0 \\
 & 200 & 5.940 & 5.199 & 0.741 & 47.8 \\
\hline
footrule & 50 & 7.583 & 6.141 & 1.442 & 55.3 \\
 & 100 & 5.854 & 5.273 & 0.581 & 40.9 \\
 & 200 & 5.503 & 4.594 & 0.909 & 64.4 \\
\hline
\end{tabular}
\end{table}

The metric graph $\Gamma_4$ replaces the 36 adjacent-swap edges by
unit intervals with shortest-path distance \citep{BakerFaber2006}.
Its vertex metric is $d_K$.
On an oriented edge $(u,v,t)$, bipartiteness gives
$d(Y,v)-d(Y,u)\in\{-1,1\}$, so
\[
 F_P(u,v,t)=(1-t)F_P(u)+tF_P(v)-t(1-t).
\]
The edge minimizer is
$t^*=\min(1,\max(0,\{1+F_P(u)-F_P(v)\}/2))$.
The minimum $329/48$ is attained at distance $1/4$ from
$\{1234,2314,3124\}$ on the six adjacent-swap edges joining this
set to $\{1324,2134,3214\}$. Their interior charts have stationary
risk and affine centred loss; full vertex support gives positive
cross-mean contrast variance. Thus \cref{prop:finite-stationary-means}
applies.

\subsection{Circular comparators}
\label{sec:supp-circle-model}
For $z_i(q)=\operatorname{wrap}(Y_i-q)\in[-\pi,\pi)$,
ordinary empirical-likelihood score inversion
\citep{KurisuOtsu2026} accepts $q$ if
\[
 -2\log\max_{\substack{p_i\ge0,\ \sum_i p_i=1\\
                         \sum_i p_i z_i(q)=0}}\prod_i(np_i)
       \le\chi^2_{1,.95}.
\]
Dual solutions are checked against the primal constraints, and
all accepted components are retained.
This fixed cutoff differs from the diverging cutoff
used for score-zero-set estimation by
\citet[Section~4.5]{KurisuOtsu2026}.

The Wald arc is centred at a global empirical
Fr\'echet mean $\widehat m_n$.
With score variance $\widehat s^2$ of $2z_i(\widehat m_n)$
(denominator $n$), its curvature, variance and radius are
\[
 \widehat H=2-4\pi\widehat f(\widehat m_n+\pi),\qquad
 \widehat V=\widehat s^2/\widehat H^2,\qquad
 r_n=\min\bigl(\pi,(\chi^2_{1,.95}\widehat V/n)^{1/2}\bigr)
 \quad(\widehat H>0).
\]
The density estimate uses a periodic von Mises kernel with
bandwidth $n^{-1/5}$. The positive-curvature intrinsic-mean
limit is given in \citet{HotzHuckemann2011}.

\subsection{Calibration and numerical optimization}
Euclidean draw suprema reduce to projection norms on the centered
sample column space. At $\delta=0$, with covariance divisor $n$,
\[
 n\bar Y_n^{\mathrm T}\wh\Sigma_n^+\bar Y_n
 \sim \frac{2n}{n-2}F_{2,n-2}.
\]
The plotted $\chi^2$ Wald comparator therefore has Gaussian coverage
$\pr(F_{2,n-2}\le(n-2)\chi^2_{2,.95}/(2n))$.
The Gaussian $F$ reference replaces the squared cutoff by
$\{2n/(n-2)\}F_{2,n-2;.95}$ and has exact 95\% coverage.
Discrete ranking enumerates all 552 ordered nontrivial pairs.
Graph and circle calibration cover all edge pairs and antipodal
chart pairs, respectively, with full-domain bounds.

For Global and Local, values at feasible comparison pairs provide
lower bounds, while full-domain bounds provide upper bounds:
\[
 \underline{S}_{n,b}^{\#}
 \le S_{n,b}^{\#}
 \le \overline{S}_{n,b}^{\#}.
\]
With $B=1999$, taking the 1900th order statistic gives the simulated
cutoff bounds and corresponding set inclusions:
\[
 \underline c_n
 \coloneqq \underline{S}_{n,(1900)}^{\#}
 \le c_n
 \le\overline{S}_{n,(1900)}^{\#}
 =:\overline c_n,
 \qquad
 \mathcal C_n(\underline c_n)
 \subseteq\mathcal C_n(c_n)
 \subseteq\mathcal C_n(\overline c_n).
\]

\renewcommand{\refname}{References for the Supplementary Material}
\putbib[ref]
\end{bibunit}
\end{document}